\documentclass[12pt]{amsart} 
\usepackage[a4paper,margin=1in,footskip=0.25in]{geometry}
\usepackage{mathrsfs}
\usepackage{amssymb}
\usepackage{amsmath}
\usepackage{amsthm}
\usepackage{mathtools}
\usepackage{nicematrix}
\usepackage{tikz-cd}
\usepackage{enumitem}
\usepackage{subfiles}
\usepackage{faktor}
\usepackage{quiver}
\usepackage{xcolor}
\usepackage[colorinlistoftodos]{todonotes}
\usepackage{xstring}
\usepackage{pdflscape}
\usepackage{subcaption}
\usepackage{float}
\usepackage{afterpage}

\usetikzlibrary{shapes.geometric,decorations.pathreplacing,decorations.markings,decorations.pathmorphing,calc,patterns,hobby}

\PassOptionsToPackage{pdfusetitle,colorlinks}{hyperref}
\usepackage{bookmark}
\hypersetup{
  linkcolor={blue!70!cyan!80!black},
  citecolor={blue!30!cyan!80!black},
  urlcolor={blue!50!cyan!80!black}
}
\usepackage{cleveref} 

\newtheorem{theorem}{Theorem}[section]
\newtheorem{lemma}[theorem]{Lemma}

\newtheorem{corollary}[theorem]{Corollary}

\theoremstyle{definition}
\newtheorem{definition}[theorem]{Definition}
\newtheorem{construction}[theorem]{Construction}

\newtheorem{warning}[theorem]{Warning}
\theoremstyle{remark}
\newtheorem{remark}[theorem]{Remark}
\newtheorem{example}[theorem]{Example}

\newcommand{\sB}{\mathcal{B}}
\newcommand{\sC}{\mathcal{C}}
\newcommand{\sD}{\mathcal{D}}
\newcommand{\sE}{\mathcal{E}}

\newcommand{\sO}{\mathcal{O}}

\newcommand{\sR}{\mathcal{R}}

\newcommand{\sT}{\mathcal{T}}

\newcommand{\bbA}{\mathbb{A}}

\newcommand{\bbN}{\mathbb{N}}

\newcommand{\bbU}{\mathbb{U}}
\newcommand{\bbZ}{\mathbb{Z}}

\newcommand{\scrB}{\mathscr{B}}
\newcommand{\scrC}{\mathscr{C}}
\newcommand{\scrD}{\mathscr{D}}

\DeclareMathOperator{\rk}{rk}

\DeclareMathOperator{\id}{id}
\DeclareMathOperator{\Hom}{Hom}

\DeclareMathOperator{\GL}{GL}

\newcommand{\GLdd}{\GL_\mathbf{d}}
\newcommand{\leqdeg}{\leq_{\textnormal{deg}}}

\newcommand{\ddim}{\textup{\textbf{dim}} \,}
\newcommand{\dd}{\textup{\textbf{d}}}
\newcommand{\hh}{\textup{\textbf{h}}}
\newcommand{\rr}{\textup{\textbf{r}}}
\newcommand{\wB}{\widetilde{B}}
\newcommand{\wm}{\widetilde{m}}
\newcommand{\wA}{\widetilde{A}}
\newcommand{\wM}{\widetilde{M}}
\newcommand{\wN}{\widetilde{N}}
\newcommand{\wlambda}{\widetilde{\lambda}}
\newcommand{\strings}{\scrC}
\newcommand{\bands}{\scrB}
\newcommand{\minbands}{\scrB_{\text{min}}}
\newcommand{\msets}{\scrD}
\newcommand{\minmsets}{\scrD_{\text{min}}}
\newcommand{\diagram}[1]{\sD(#1)}

\let\mod\relax
\DeclareMathOperator{\mod}{mod}

\let\amsamp=&
\newcommand{\pmatquiver}[1]{%
  \left(\begin{smallmatrix}#1\end{smallmatrix}\right)%
}

\newcommand{\stackquiver}[1]{%
  \begin{matrix}#1\end{matrix}%
}

\newcommand{\pmat}[1]{%
\begin{pNiceMatrix}
[small]
#1
\end{pNiceMatrix}
}

\newcommand{\biblio}{\bibliographystyle{amsalpha}\bibliography{literature}}

\pgfdeclarelayer{background}
\pgfsetlayers{background,main}

\tikzset{
  on each segment/.style={
    decorate,
    decoration={
      show path construction,
      moveto code={},
      lineto code={
        \path [#1]
        (\tikzinputsegmentfirst) -- (\tikzinputsegmentlast);
      },
      curveto code={
        \path [#1] (\tikzinputsegmentfirst)
        .. controls
        (\tikzinputsegmentsupporta) and (\tikzinputsegmentsupportb)
        ..
        (\tikzinputsegmentlast);
      },
      closepath code={
        \path [#1]
        (\tikzinputsegmentfirst) -- (\tikzinputsegmentlast);
      },
    },
  },
  mid arrow/.style={postaction={decorate,decoration={
        markings,
        mark=at position .5 with {\arrow[#1]{>}}
      }}},
}

\tikzset{snake it/.style={decorate, decoration=snake}}

\def\centerarc[#1](#2)(#3:#4:#5)
{ \draw[#1] ($(#2)+({#5*cos(#3)},{#5*sin(#3)})$) arc (#3:#4:#5) }

\newcommand*\circled[1]{\tikz[baseline=(char.base)]{
            \node[shape=circle,draw,inner sep=2pt] (char) {#1};}}

\newcommand*{\drawStringOrBand}[3]{
    \begin{tikzpicture}[scale = 0.15]
        \foreach \x [count=\xcount] in {#2}{
            \IfBeginWith{\x}{a}{\draw (\xcount, 2) -- (\xcount + 1, 0);}{
                \IfBeginWith{\x}{b}{\draw (\xcount, 0) -- (\xcount + 1, 2);}{
                    \IfBeginWith{\x}{c}{\draw (\xcount, -2) -- (\xcount + 1, 0);}{
                        \IfBeginWith{\x}{d}{\draw (\xcount, 0) -- (\xcount + 1, -2);}{}
                    }
                }
            }
        }
        \draw[dotted] (1,0) -- (#3 + 1,0);
        \path (1,-2) -- (1,2); 
        \IfBeginWith{#1}{1}{
            \IfBeginWith{#2}{a}{
                \node[circle,fill=black,inner sep=0pt,minimum size=3pt] at (1, 2) {};
                \node[circle,fill=black,inner sep=0pt,minimum size=3pt] at (#3 + 1, 2) {};
            }{
                \IfBeginWith{#2}{c}{
                    \node[circle,fill=black,inner sep=0pt,minimum size=3pt] at (1, -2) {};
                    \node[circle,fill=black,inner sep=0pt,minimum size=3pt] at (#3 + 1, -2) {};
                }{
                    \node[circle,fill=black,inner sep=0pt,minimum size=3pt] at (1, 0) {};
                    \node[circle,fill=black,inner sep=0pt,minimum size=3pt] at (#3 +1 , 0) {};
                }
            }
        }{}
   \end{tikzpicture}
}

\newcommand*{\drawStringOrBandKronecker}[3]{
    \begin{tikzpicture}[scale = 0.15]
        \foreach \x [count=\xcount] in {#2}{
            \IfBeginWith{\x}{a}{\draw (\xcount, 2) -- (\xcount + 1, 0);}{
                \IfBeginWith{\x}{b}{\draw (\xcount, 0) -- (\xcount + 1, 2);}{
                }
            }
        }
        \IfBeginWith{#1}{1}{
            \IfBeginWith{#2}{a}{
                \node[circle,fill=black,inner sep=0pt,minimum size=3pt] at (1, 2) {};
                \node[circle,fill=black,inner sep=0pt,minimum size=3pt] at (#3 + 1, 2) {};
            }{
                    \node[circle,fill=black,inner sep=0pt,minimum size=3pt] at (1, 0) {};
                    \node[circle,fill=black,inner sep=0pt,minimum size=3pt] at (#3 +1 , 0) {};
            }
        }{}
   \end{tikzpicture}
}

\newcommand{\leqdegfam}{\leq_{\text{deg}}}
\newcommand{\df}[1]{{\color{blue!50!teal}\textit{#1}}} 

\begin{document}

\renewcommand{\biblio}{}

\title{Degenerations of families of bands and strings for gentle algebras}

\author{Judith Marquardt}
\address[Judith Marquardt]{Univ. Grenoble Alpes, CNRS, IF, 38000 Grenoble, France, and Université Paris-Saclay, UVSQ, CNRS, Laboratoire de Mathématiques de Versailles, 78000, Versailles, France}
\email{\href{mailto:judith.marquardt@univ-grenoble-alpes.fr}{judith.marquardt@univ-grenoble-alpes.fr}}

\keywords{Representation Theory, Gentle Algebras, Degeneration}
\thanks{}

\begin{abstract}
Let $A$ be a gentle algebra. For every collection of string and band diagrams, we consider the constructible subset of the variety of representations containing all modules with this underlying diagram. We study degenerations of such sets. We show that these sets are defined by vectors of integers which we call $h$-vectors and which are related to a restricted version of the $\hom$-order. We provide combinatorial criteria for the existence of a degeneration, involving the removal of an arrow or the resolving of a type of configuration called ``reaching''.
\end{abstract}

\maketitle

\tableofcontents

\biblio

\section{Introduction}
\label{sect:introduction}

Degenerations of orbits in a module variety are well studied, see for example \cite{Zwara2000}. Gentle algebras were introduced by \cite{Assem_Skowronski_intro_gentle}. They are tame and their modules are completely classified \cite{Wald_Waschbüch_indecomposable_gentle,ButlerRingel87}. Recently, they appeared in the context of surface combinatorics, see for example \cite{HKK,OPS,PPP_kissing_locally,B-CS}. Their module varieties have also been the object of recent results: \cite{GLS_gentle_2021} systematically studied their irreducible components and generically $\tau$-reduced components. There are also related results in the direction of special biserial algebras by \cite{CCKW}.

The isomorphism classes of indecomposable modules of a gentle algebra are governed by the combinatorics of strings and bands of the quiver. To every string we associate an orbit of a module, to every band a family of orbits. In this article, for every diagram $\sD$ consisting of strings and bands, we consider the set $\sO_\sD$ of all modules whose underlying diagram is $\sD$, see \Cref{def:families_min_bands}. We call this the family of $\sD$. This ``combinatorial'' definition has a geometrical sense. As shown by \cite{GLS_gentle_2021}, the indecomposable blocks of the module variety of an algebra are the closure of $\sO_\sD$ where $\sD$ is either a single string or a single band. Motivated by these results, we initiate a systematic study of closures of the sets $\sO_\sD$.

We introduce a degeneration order $\sD \leqdegfam \sD'$ on the diagrams, defined by $\sO_{\sD'} \subseteq \overline{\sO_\sD}$, see \Cref{def:degeneration_of_strings_and_bands}. The main objective of this article is to study this order.

In the classical study of degenerations of orbits of modules, several orders are used. Recall that the degeneration order is defined by $M \leqdeg M'$ if $\sO_{M'} \subseteq \overline{\sO_M}$. It implies a weaker order given by assigning an infinite vector of integers to each module $M$ and comparing them pointwise. This is called the $\hom$-order, defined by $M \leq_{\hom} M'$ if $\hom(Z,M) \leq \hom(Z,M')$ for all modules $Z$. It is known that $M\leqdeg M'$ implies $M\leq_{\hom} M'$ while the inverse is not always true, see for example \cite{CB}.

We adapt the $\hom$-order to the study of the families $\sO_\sD$. We do this by restricting the $\hom$-order to comparisons involving only string modules. Precisely, we define the $h$-order by ${M \leq_h M'}$ if ${\hom (Z,M) \leq \hom(Z,M')}$ for all string modules $Z$. We call ${h(M) := (\hom(Z,M))_{Z \text{ string module}}}$ the $h$-vector of $M$.

Our main results are as follows.
\begin{itemize}
    \item The sets $\sO_\sD$ are precisely the subsets of the module variety on which the $h$-vector is constant, see \Cref{cor:families_h}.
    \item The degeneration order $\sD \leqdegfam \sD'$ implies the $h$-order $\sD \leq_h \sD'$, see \Cref{lem:two_partial_orders}. We do not know whether the converse holds.
\end{itemize}
Next, we provide combinatorial criteria which imply degenerations.
\begin{itemize}
    \item The deletion of an arrow in a string or band, see \Cref{lem:deleting_arrow}.
    \item The resolution of a ``reaching'', a configuration of strings and band described in \cite{PPP_kissing}. Resolving is an explicit combinatorial operation that implies a degeneration, see \Cref{thm::deg_and_intersection}.
\end{itemize}
In all examples that we have computed, cover relations of the degeneration order always fall in the above two cases. We do not know whether this is true in general.

This paper is organised as follows. We introduce some notation and basic concepts of module varieties and gentle algebras in \Cref{sect:preliminaries}. We define the sets $\sO_\sD$ in \Cref{sect:families} and study the $h$-vectors. Then we introduce degeneration of strings and bands and study its existence in \Cref{sect:degeneration}. In \Cref{sect:examples} we discuss some examples. Finally, \Cref{sect:outlook} contains open questions and possible directions for future research.

\biblio

\section{Preliminaries}
\label{sect:preliminaries}

We recall some necessary notation. For more details, we refer to \cite{ASS} for an introduction into quiver representations, to \cite{PPP_kissing} for gentle algebras and the formalism of strings and bands and to \cite{Zwara_Survey} for an overview of module varieties. We work over an algebraically closed field $K$.

\subsection{Path algebras and representations}

A \df{quiver} is a directed graph given by the data $Q = (Q_0, Q_1, s, t)$ where $Q_0$ denotes the set of vertices and $Q_1$ the set of arrows. The maps $s,t \colon Q_1 \rightarrow Q_0$ assign the source resp. target of each arrow. We assume quivers to be connected and finite and set $n$ to be the number of vertices.
We compose arrows as we would compose morphisms. A \df{path} $p = a_1 \ldots a_m$ is a composition of compatible arrows with $t(a_{i+1}) = s(a_i)$ for $1 \leq i < m$. The \df{length} of $p$ is $m$ and denoted by $\ell(p)$. We extend the maps $s$ and $t$ to paths by setting $s(p) = s(a_m)$ and $t(p) = t(a_1)$. Moreover, for every vertex $i \in Q_0$ we define a 0-length path $e_i$ with $s(e_i) = t(e_i) = i$.

Recall that a quiver $Q$ gives rise to the \df{path algebra} $KQ$ which is the $K$-algebra with basis given by all paths of $Q$ and an algebra structure induced by the composition of paths. Let $\mathfrak{m}$ be the ideal generated by the arrows of $Q$. We say that an ideal $I$ of $KQ$ is \df{admissible} if there exists an integer $m \geq 2$ such that $\mathfrak{m}^m \subseteq I \subseteq \mathfrak{m}^2$. A \df{basic} algebra is an algebra isomorphic to $A = KQ/I$ where $I$ is admissible. Since we work over an algebraically closed field, every finite-dimensional algebra is Morita equivalent to a basic algebra. A \df{representation} of $A$ is a tuple $(M_i,M_a)_{i \in Q_0,\, a \in Q_1}$ where $M_i$ is a finite-dimensional $K$-vector space and $M_a \colon M_{s(a)}\rightarrow M_{t(a)}$ is a $K$-linear morphism such that the relations of $I$ are satisfied. The \df{dimension vector} of $M$ is $\ddim(M) = (\dim M_1, \ldots, \dim M_n )$. We identify representations of $A$ with finite-dimensional left modules of $A$ in the usual way.

\subsection{Gentle algebras}  

\begin{definition}
    \label[definition]{def_gentle_alg}
    A \df{gentle algebra} is a basic algebra $A = {KQ}/{I}$ such that
    \begin{itemize}
        \item every vertex of $Q$ admits at most two incoming arrows and two outgoing arrows,
        \item the ideal $I$ is generated by paths of length two,
        \item for arrows $a$, $b$, $c$ with $a \neq b$ and $t(a) = t(b) = s(c)$, exactly one of the paths $ca$, $cb$ is in $I$,
        \item dually, for arrows $a$, $b$, $c$ with $a \neq b$ and $s(a) = s(b) = t(c)$, exactly one of the paths $ac$, $bc$ is in $I$.
    \end{itemize}
\end{definition}

\begin{example}
\label[example]{easy_example_gentle_algebra}
Consider the following quiver $Q$. The algebra $A = {KQ}/\langle{ba}\rangle$ is gentle. The relation $ba$ is indicated by the dashed line. Note that the path algebra without relations is not gentle.
\[\begin{tikzcd}[sep = small]
	&& 3 \\
	Q=~1 & 2 \\
	&& 4
	\arrow["c"', from=2-2, to=3-3]
	\arrow[""{name=0, anchor=center, inner sep=0}, "b"', from=2-2, to=1-3]
	\arrow[""{name=1, anchor=center, inner sep=0}, "a"', from=2-1, to=2-2]
	\arrow[curve={height=-12pt}, shorten <=7pt, shorten >=7pt, dashed, no head, from=1, to=0]
\end{tikzcd}\]
\end{example}

For the remainder of this work, let $A = KQ/I$ be a gentle algebra.

\subsection{Strings and bands}

The modules of gentle algebras have been intensely studied in \cite{ButlerRingel87, Wald_Waschbüch_indecomposable_gentle}. They are closely related to the notion of strings and bands.

We define $Q_1^-:=\left\{a^- \,|\, a \in Q_1 \right\}$ where $a^-$ is the formal inverse arrow of $a$ with $s(a^-) = t(a)$ and $t(a^-) = s(a)$. Moreover, we set $(a^-)^- := a$. This allows us to define the \df{double quiver} $Q^{\pm} = ( Q_0,  Q_1 \cup Q_1^-, s, t)$. For a path $p = c_1 \dots c_m$ in $Q^{\pm}$, we define its inverse $p^-:= c_m^- \dots c_1^-$. Let $S_I$ be the generating set of $I$ consisting of length two paths. We define $I^{\pm} := \langle S_I \cup \{ p^- \,|\, p \in S_I \} \rangle$ which is an ideal of $KQ^{\pm}$.

\begin{definition}
\label[definition]{def_string}
A \df{string} $C$ is a path $c_1 \dots c_m$ in the double quiver such that
\begin{enumerate}[label=(\roman*)]
    \item no arrow is followed directly by its inverse, i.e. $aa^-$ and $a^-a$ are not allowed to be a part of the path $\forall a \in Q_1$,
    \item $C \notin I^\pm$.
\end{enumerate}
We define the \df{length} of $C$ as the length of the path and denote it by $\ell(C)$. We handle $s$ and $t$ the same way. Note that the inverse of $C$, $C^-$, is also a string. 
\end{definition}

We consider strings up to the inverse action. So a string is an element of $\strings = \{\text{strings}\}/_{C\sim C^-}$ and we choose a representative $C$ to depict a string $\sC$.

For an arrow $a\in Q_1$, $ i \xrightarrow[]{a} j$, we depict its appearance within a string as

\begin{center}
    \begin{tabular}{ccc}
        $\begin{tikzcd}[sep=small]
	       & i \\
	       j
	       \arrow["a"', from=1-2, to=2-1]
    \end{tikzcd}$ 
    & \qquad \qquad &
    $\begin{tikzcd}[sep=small]
	       i  \\
	       & j
	       \arrow["a", from=1-1, to=2-2]
    \end{tikzcd}$ 
    \\
    for $a$ & & for $a^-$.
    \end{tabular}
\end{center}

\begin{example}
\label[example]{ex:string}
    We consider the algebra from Example \ref{easy_example_gentle_algebra}. Then $ca$ and $cb^-$ are strings and depicted in the following way
    \[\begin{tikzcd}[sep=small]
	&& 1 \\
	& 2 \\
	4
	\arrow["a"', from=1-3, to=2-2]
	\arrow["c"', from=2-2, to=3-1]
    \end{tikzcd}
    \qquad \qquad
    \begin{tikzcd}[sep=small]
	& 2 \\
	4 && 3.
	\arrow["c"', from=1-2, to=2-1]
	\arrow["b", from=1-2, to=2-3]
\end{tikzcd}
    \]
    Note that $ca$ can also be represented by $a^-c^-$ and $cb^-$ by $bc^-$.
\end{example}

\begin{remark}
    Since for a representative $C$, the source lies on the right and the target on the left, we will often use these directions for orientation within the string. 
\end{remark}

\begin{definition}
\label[definition]{def:band}
A \df{band} is a string $B$ of length at least 2 such that $B^2$ is also a string. We call a band $B$ \df{minimal} if $B$ is minimal as a cycle, i.e. there exists no string $C$ and $s \geq 2$ such that $C^s = B$. A band that is not minimal will be called \df{non-minimal}.
\end{definition}

\begin{warning}
    In most of the literature, the term band is used for what we will call minimal band. In the context of this work, it is often easier to include non-minimal bands. We will usually denote them as $B^t$ where $B$ is a minimal band and $t>1$. Note that at the end of \Cref{sect:families,sect:degeneration}, we will restrict to minimal bands.
\end{warning}

We consider bands up to inverse action and cyclic permutation. As for strings, we define the set $\bands = \{\text{bands}\}/_{\sim}$. Then we consider a band $\sB$ to be an element of $\bands$, represented by a specific orientation $B$. We similarly define the set of minimal bands $\minbands \subseteq \bands$.

\begin{example}
\label[example]{ex:band}
    Let $Q$ be the Kronecker quiver 
    $\begin{tikzcd}
	   1 & 2.
	   \arrow["a", shift left, from=1-1, to=1-2]
	   \arrow["b"', shift right, from=1-1, to=1-2]
    \end{tikzcd}$
    The path algebra $KQ$ admits a unique minimal band $a^-b$. We depict bands by circling the first and last vertex:
    \begin{equation*}
        \begin{tikzcd}[sep=small]
	        \circled{1} && \circled{1}. \\
        	& 2
	        \arrow["a", from=1-1, to=2-2]
        	\arrow["b"', from=1-3, to=2-2]
        \end{tikzcd}
    \end{equation*}
    This allows us to visually separate strings and bands. Note that the same band is also represented by $ba^-$, $ab^-$ and $b^-a$.
\end{example}

\begin{definition}
\label[definition]{def:infband}
    Let $B$ be a band. Then we associate to it a string of infinite length, the \df{infinitely unraveled band} $B^\infty$, which consists of $\bbZ$ copies of $B$. 
\end{definition}

\begin{example}
\label[example]{ex:unraveled_band}
    Let $B$ be the band from \Cref{ex:band}. Then $B^\infty$ is represented by
    \[\begin{tikzcd}[sep=small]
    	& 1 && 1 && 1 && 1 && 1 && 1 \\
    	{} && 2 && 2 && 2 && 2 && 2. && {}
    	\arrow["\cdots"', draw=none, from=1-2, to=2-1]
    	\arrow["a", from=1-2, to=2-3]
    	\arrow["b"', from=1-4, to=2-3]
    	\arrow["a", from=1-4, to=2-5]
    	\arrow["b"', from=1-6, to=2-5]
    	\arrow["a", from=1-6, to=2-7]
    	\arrow["b"', from=1-8, to=2-7]
    	\arrow["a", from=1-8, to=2-9]
    	\arrow["b"', from=1-10, to=2-9]
    	\arrow["a", from=1-10, to=2-11]
    	\arrow["b"', from=1-12, to=2-11]
    	\arrow["\cdots", draw=none, from=1-12, to=2-13]
    \end{tikzcd}\]
\end{example}

We will need some additional terminology of substrings of strings and bands. We follow \cite{PPP_kissing}.

\begin{definition}
\label[definition]{def:substring}
    A \df{substring} of a string $C = c_1 \ldots c_m$ is a subword of consecutive arrows of the string. Every substring has the form $C[i,j] = c_i \ldots c_{j-1}$ where $1 \leq i \leq j \leq m+1$. The $0$-length substrings are $C[i,i] = e_{t(a_i)}$ for all $1 \leq i \leq m$ and $C[m+1,m+1] = s(a_m)$. Note that the position of the substring in $C$ is a part of its data. We denote $C[1,i]$ by $C_{\leq i}$ and $C[j,m+1]$ by $C_{\geq j}$.
\end{definition}

\begin{example}
\label[example]{ex:substring}
    We return to the Kronecker quiver in \Cref{ex:band} and consider the string $C = ba^-b$. \Cref{fig:substrings} depicts the notation for substrings. The bottom row shows the values for $i$ and $j$ for a substring between these vertices. Note that we number vertices from tail to source. The string $b$ appears twice as a substring, once as $C_{\leq 2} = C[1,2]$ and once as $C_{\geq 3} = C[3,4]$. These do not coincide as substrings.
    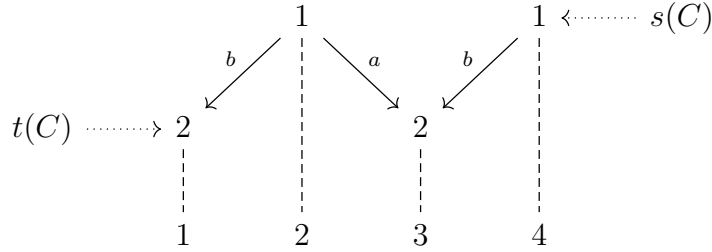
\begin{figure}[H]
        \centering
        \[\begin{tikzcd}
        	&& 1 && 1 & {s(C)} \\
        	{t(C)} & 2 && 2 \\
        	& 1 & 2 & 3 & 4
        	\arrow["b"', from=1-3, to=2-2]
        	\arrow["a", from=1-3, to=2-4]
        	\arrow[dashed, no head, from=1-3, to=3-3]
        	\arrow["b"', from=1-5, to=2-4]
        	\arrow[dashed, no head, from=1-5, to=3-5]
        	\arrow[dotted, from=1-6, to=1-5]
        	\arrow[dotted, from=2-1, to=2-2]
        	\arrow[dashed, no head, from=2-2, to=3-2]
        	\arrow[dashed, no head, from=2-4, to=3-4]
        \end{tikzcd}\]
        \caption{Notation of substrings.}
        \label{fig:substrings}
    \end{figure}
\end{example}

\begin{definition}
\label[definition]{def:substring_unraveled_band}
    Let $B$ be a band. A \df{substring} of the infinitely unraveled band $B^\infty$ is a finite substring of $B^\infty$ as in \Cref{def:substring}. Substrings of $B^\infty$ are considered up to shifts by $\ell(B)$, i.e. $B^\infty[i,j] = B^\infty[i+\ell(B), j+\ell(B)]$.
\end{definition}

\begin{example}
\label[example]{ex:substring_unraveled_band}
    Following \Cref{ex:band,,ex:unraveled_band}, $B^\infty$ contains infinitely many arrows $b$ but they all give rise to the same substring of $B^\infty$. 
\end{example}

\begin{definition}
\label[definition]{def:substring_band}
    Let $B$ be a band and $\rho$ a substring of $B^\infty$. We let $\rho^\text{red}$ be the shortest substring at the beginning of $\rho$ such that the remaining part of $\rho$ is equal to some power of $B$ as bands. We call $\rho_\text{red}$ the \df{substring of $B$ associated to $\rho$}. Note that there is always a representative of $B$ so that we can write $\rho_\text{red}$ as $B[i,j]$ with $1 \leq i \leq j \leq \ell(B)+1$.
\end{definition}

\begin{example}
\label[example]{ex:substring_band}
    We consider the Kronecker quiver and $B = a^-b$ as before, see \Cref{ex:band}. Then $ba^-b$ is a substring of $B^\infty$ and the associated substring of $B$ is $b = B[2,3] = B_{\geq 2}$.
\end{example}

\begin{definition}
\label[definition]{def:oriented_substrings}
    Let $C$ be a string with substrings $C[i,j]$ and $C[i',j']$ of length at least 1 that coincide as strings. We say that the substrings are \df{oriented in the same direction} in $C$ if the representative $C[i,j]$ is the same as the representative $C[i',j']$. Otherwise, we call them \df{oriented in opposite directions}. We extend this definition to substrings of unraveled bands.
\end{definition}

\begin{remark}
    Note that two substrings are either oriented in the same direction or in opposite directions as strings cannot be palindromes.
\end{remark}

\begin{example}
    We consider \Cref{ex:substring}. The substrings $C_{\leq 2}$ and $C_{\geq 3}$ are oriented in the same direction.
\end{example}

\begin{definition}
\label[definition]{def:top_string}
    Let $C$ be a string. A substring $C[i,j]$ is \df{on top of} $C$ if the following two conditions are fulfilled.
    \begin{itemize}
        \item There is no arrow left of the substring or the next arrow to the left of the substring lies in $Q_1$, i.e. $i = 1$ or $a_{i-1} \in Q_1$.
        \item There is no arrow right of the substring or the next arrow to the right of the substring lies in $Q_1^-$, i.e. $j= m+1$ or $a_j \in Q_1^-$.
    \end{itemize} 
    We denote the set of substrings on top of $C$ as $\Sigma_\text{top}(C)$. For a string $E$, we denote the subset of substrings on top of $C$ of the form $E$ by $\Sigma_\text{top}^E(C)$.
\end{definition}

\begin{example}
\label[example]{ex:top_string}
    We reconsider \Cref{ex:substring}. The left copy of $b$, $C_{\leq 2}$ is on top of $C$. However, the right copy, $C_{\geq 3}$, is not. An example that has neighbouring arrows on both sides is $e_1 = C[2,2]$ which is on top of $C$. Note that $C$ itself also is on top of $C$.
\end{example}

\begin{definition}
\label[definition]{def:top_band}
    Let $B$ be a band. A substring of $B^\infty$ is \df{on top of} $B^\infty$ if the next arrow to the left of it is in $Q_1$ and the next arrow to the right of it is in $Q_1^-$.
    We denote the set of substrings on top of $B^\infty$ as $\Sigma_\text{top}(B^\infty)$. For a string $E$, we denote the subset of substrings on top of $B^\infty$ of the form $E$ by $\Sigma_\text{top}^E(B^\infty)$.
\end{definition}

\begin{example}
\label[example]{ex:top_band}
    We consider \Cref{ex:unraveled_band} and the substring $a^-b$ of $B^\infty$. It is on top of $B^\infty$.
\end{example}

Dually, we define substrings at the bottom of a string or unraveled band.

\begin{definition}
\label[definition]{def:bot_string}
    Let $C$ be a string. A substring $C[i,j]$ is \df{at the bottom of} $C$ if the following two conditions are fulfilled.
    \begin{itemize}
        \item There is no arrow left of the substring or the next arrow to the left of the substring lies in $Q_1^-$, i.e. $i = 1$ or $a_{i-1} \in Q_1^-$.
        \item There is no arrow right of the substring or the next arrow to the right of the substring lies in $Q_1$, i.e. $j= m+1$ or $a_j \in Q_1$.
    \end{itemize} 
    We denote the set of substrings at the bottom of $C$ as $\Sigma_\text{bot}(C)$. For a string $E$, we denote the subset of substrings at the bottom of $C$ of the form $E$ by $\Sigma_\text{bot}^E(C)$.
\end{definition}

\begin{example}
\label[example]{ex:bot_string}
    We reconsider \Cref{ex:substring}. The right copy of $b$, $C_{\geq 3}$ is at the bottom of $C$. Note that the substring $a = C[2,3]$ is neither on top nor at the bottom of $C$.
\end{example}

\begin{definition}
\label[definition]{def:bot_band}
    Let $B$ be a band. A substring of $B^\infty$ is \df{at the bottom of} $B^\infty$ if the next arrow to the left of it is in $Q_1^-$ and the next arrow to the right of it is in $Q_1$.
    We denote the set of substrings at the bottom of $B^\infty$ as $\Sigma_\text{bot}(B^\infty)$. For a string $E$, we denote the subset of substrings at the bottom of $B^\infty$ of the form $E$ by $\Sigma_\text{bot}^E(B^\infty)$.
\end{definition}

\begin{example}
\label[example]{ex:bot_band}
    We consider \Cref{ex:unraveled_band}. The substring $ba^-$ is at the bottom of $B^\infty$.
\end{example}

\begin{remark}
\label[remark]{rmk:band_rota_top_bot}
    Let $B$ be a band. Since $A$ is finite-dimensional, $B$ must contain arrows in $Q_1$ and in $Q_1^-$. In particular, there is a representative of $B$ such that, considered as a string, it is a substring at the top of $B^\infty$. Dually, there is a representative that gives a substring at the bottom of $B^\infty$. This can be seen in \Cref{ex:top_band,,ex:bot_band}.
\end{remark}

We will later use substrings to distinguish certain strings and bands. For this, we observe the following.

\begin{lemma}
\label[lemma]{lem:bands_substring_length}
    Let $B$ and $B'$ be two distinct minimal bands. Then there exists an integer $m$ such that $\forall n \geq m$, for every representative $\wB$ of $B$, $(\wB)^n \notin \Sigma_{\text{bot}}((B')^\infty)$.
\end{lemma}

\begin{proof}
    First, let us show that it is enough to prove that there exists an integer $m$ such that no representative of $B^m$ is a substring at the bottom of $(B')^\infty$. Let us assume there exists an $m$ and a representative $\wB$ of $B$ such that $(\wB)^m \in \Sigma_{\text{bot}}((B')^\infty)$. We let $B_\text{bot}$ be a representative of $B$ such that it is at the bottom of $B^\infty$, see \Cref{rmk:band_rota_top_bot}. Then for every $n < m$, $B_\text{bot}^n$ is a substring at the bottom of $(B')^\infty$. So finding an $m$ where no representative of $B^m$ is such a substring will imply the same for all higher integers.
    
    Now let $\ell = \ell(B)$, $\ell' = \ell(B')$. If there exists a representative $\wB$ of $B$ such that $(\wB)^{\ell'} \in \Sigma_{\text{bot}}((B')^\infty)$ then $(\wB)^{\ell'} = (B')^\ell$. By minimality of the bands, this implies $B = B'$.
\end{proof}

\subsection{String and band modules}

Every string gives rise to an indecomposable $A$-module. Every (minimal) band gives rise to families of (indecomposable) $A$-modules. We follow the construction of \cite{Wald_Waschbüch_indecomposable_gentle}.

\begin{construction}[String modules]
\label[construction]{constr::strings}
    Let $C = c_1 \dots c_m$ be a string. We will construct its associated \df{string module} denoted by $M(C)$. We let $A_{m+1}$ be the path algebra of the quiver of Dynkin type $A_{m+1}$ 
    \[\begin{tikzcd}
	1 & 2 & \ldots & m & {m+1}
	\arrow["a_1",no head, from=1-1, to=1-2]
	\arrow[no head, from=1-2, to=1-3]
	\arrow[no head, from=1-3, to=1-4]
	\arrow["a_m", no head, from=1-4, to=1-5]
\end{tikzcd}\]
    whose orientation of arrows is induced by $C$ in the following way. For all $1 \leq i \leq m$ the  arrow $a_i$ between the vertices $i$ and $i+1$ is oriented to the right, $i \xlongrightarrow{a_i} i+1$, if $c_i \in Q_1^-$. It is oriented to the left, $i \xlongleftarrow{a_i} i+1$, if $c_i \in Q_1$. There is a quiver morphism 
    \begin{align*}
        f \colon A_{m+1} &\rightarrow Q \\
        a_i &\mapsto 
        \left\{
        \begin{aligned}
            c_i &\quad \textnormal{ if } c_i \in Q_1, \\
            c_i^- &\quad \textnormal{ if } c_i \in Q_1^-. 
        \end{aligned}
        \right.
    \end{align*}
     This morphism induces a functor $F_f: \mod(A_{m+1}) \rightarrow \mod(A)$. We consider the $A_{m+1}$-module $M = (K, \id)_{i,a_j}$ and define $M(C) := F_f(M)$. In particular, $M(C)$ is an $(m+1)$- dimensional module.
\end{construction}

\begin{example}
\label[example]{ex:string_module}
    We consider the Kronecker quiver, see \Cref{ex:band}, and the string $C = a$. Then $M(a)$ is the representation
    \[\begin{tikzcd}
    	K & K.
    	\arrow["1", shift left=1, from=1-1, to=1-2]
    	\arrow["0"', shift right=1, from=1-1, to=1-2]
    \end{tikzcd}\]
\end{example}

\begin{construction}[Band modules]
\label[construction]{constr::bands}
    For a band $B = b_1 \ldots b_m$, $\lambda \in K^*$ and $q \in \bbN_{>0}$, the \df{band module} $M(B, \lambda, q)$ is constructed as follows. We let $\widetilde{A}_{m-1}$ be the path algebra of the Euclidean graph $\widetilde{A}_{m-1}$
    \[\begin{tikzcd}
	1 & 2 & \ldots & {m-1} & m
	\arrow["a_1", no head, from=1-1, to=1-2]
	\arrow["a_m"', curve={height=30pt}, no head, from=1-1, to=1-5]
	\arrow[no head, from=1-2, to=1-3]
	\arrow[no head, from=1-3, to=1-4]
	\arrow["a_{m-1}", no head, from=1-4, to=1-5]
\end{tikzcd}\]
    whose orientation of arrows is induced by $B$ as before. Again, there is a quiver morphism $f \colon \widetilde{A}_{m-1} \rightarrow Q$ inducing a functor $F_f: \mod(\widetilde{A}_{m-1}) \rightarrow \mod(A)$. We consider the $\widetilde{A}_{m-1}$-module $M_{q,\lambda} = (K^q, f_j)_{i,a_j}$ where 
    \begin{equation*}
        f_j = \left\{
        \begin{aligned}
            \id &\quad \textnormal{ if } j \neq m, \\
            J(q,\lambda) &\quad \textnormal{ if } j = m.
        \end{aligned}
        \right.
    \end{equation*}
    Here, $J(q,\lambda)$ denotes the Jordan block of size $q$ with eigenvalue $\lambda$. In particular, if $q = 1$, $f_m = \lambda$. Now we define $M(B, \lambda, q) := F_f(M_{q,\lambda})$. It is an $mq$-dimensional module. We call $q$ the \df{quasi-length} of $M(B, \lambda, q)$.
\end{construction}

\begin{example}
\label[example]{ex:band_module}
    Recall \Cref{ex:band}. For $\lambda \in K$, the module $M(b^-a, \lambda, 1)$ corresponds to the representation
    \[\begin{tikzcd}
    	K & K.
    	\arrow["\lambda", shift left=1, from=1-1, to=1-2]
    	\arrow["1"', shift right=1, from=1-1, to=1-2]
    \end{tikzcd}\]
    The module $M(b^-a,\lambda,2)$ is given by
    \[\begin{tikzcd}
    	{K^2} & {K^2}.
    	\arrow["{\pmatquiver{\lambda \amsamp 1 \\ \amsamp \lambda}}", shift left=1, from=1-1, to=1-2]
    	\arrow["{\pmatquiver{1 \amsamp  \\ \amsamp 1}}"', shift right=1, from=1-1, to=1-2]
    \end{tikzcd}\]
\end{example}

Note that the above construction is also applicable to non-minimal bands. The resulting modules are decomposable into minimal band modules.

\begin{example}
\label[example]{ex:non_min_band_module}
    We reconsider \Cref{ex:band_module}. Let $B$ be the minimal band $b^-a$. Then $B^2 = b^-ab^-a$ is a non-minimal band. For $\lambda \in K$, the module $M((b^-a)^2, \lambda, 1)$ corresponds to the representation
    \[\begin{tikzcd}
    	{K^2} & {K^2}.
    	\arrow["{\pmatquiver{ \amsamp \lambda \\ 1 \amsamp }}", shift left=1, from=1-1, to=1-2]
    	\arrow["{\pmatquiver{1 \amsamp  \\ \amsamp 1}}"', shift right=1, from=1-1, to=1-2]
    \end{tikzcd}\]
    This module is isomorphic to the direct sum $M(B, \sqrt{\lambda},1) \oplus M(B, -\sqrt{\lambda},1)$. Moreover, it is not isomorphic to $M(B, \lambda, 2)$.
\end{example}

\begin{theorem}[\cite{ButlerRingel87, Wald_Waschbüch_indecomposable_gentle}]
\label{thm::ind_mod_of_gen_alg} 
    Let $A = KQ / I$ be a gentle algebra. The modules $M(C)$ and $M(B,\lambda,q)$ with $C \in \strings$ , $B \in \minbands$, $\lambda \in K^*$ and $q \geq 1$ form a complete set of pairwise non-isomorphic representatives of isomorphism classes of indecomposable modules in $\mod(A)$.
\end{theorem}

\begin{definition}
\label[definition]{def:diagram}
    Let $M$ be a module. Then there is a unique decomposition of $M$ into indecomposable modules of the form
    \begin{equation*}
        M \cong M(C_1) \oplus \ldots \oplus M(C_s) \oplus M(B_1,\lambda_1,q_1) \oplus \ldots \oplus M(B_t,\lambda_t,q_t).
    \end{equation*}
    We call the \df{diagram} of $M$ the multi-set of strings and minimal bands 
    \begin{equation*}
        \diagram{M} = \{C_1, \ldots, C_s, B_1^{\times q_1}, \ldots, B_t^{\times q_t} \}
    \end{equation*}
    where we use the notation $B_i^{\times q_i}$ to mean that $q_i$ copies of $B_i$ belong to the multi-set. 
\end{definition}

\begin{example}
\label[example]{ex:diagram}
    We consider the Kronecker quiver, $\lambda \in K^*$ and the module $M$ of the form
    \[\begin{tikzcd}
    	K^2 & K^2.
    	\arrow["{\left(\begin{smallmatrix} \lambda \amsamp \\ \amsamp 1 \end{smallmatrix}\right)}", shift left=1, from=1-1, to=1-2]
    	\arrow["{\left(\begin{smallmatrix} 1 \amsamp \\ \amsamp 0 \end{smallmatrix}\right)}"', shift right = 1, from=1-1, to=1-2]
    \end{tikzcd}\]
    Then $M = M(b^-a,\lambda, 1) \oplus M(a)$, see \Cref{ex:band_module,,ex:string_module}. So its diagram is $\diagram{M} = \{C,B\}$ with $C = a$, $B = b^-a$.
\end{example}

\begin{remark}
\label[remark]{rmk:notation_multi_sets}
    We use three different notations for multi-sets depending on the context. The first is the short form $\{a_i^{\times q_i}\}_i$ where the $a_i$ are pairwise different elements and the $q_i$ denote their multiplicities. The second is the long form, $\{a_i\}_i$, written as a set where elements are not necessarily pairwise distinct. The third form is a mix of the first two, as seen in \Cref{def:diagram}.
    
    We let $\msets$, resp.$\minmsets$, be the set of multi-sets of bands and strings, resp. minimal bands and strings. For one-element multi-sets, we may drop the parenthesis, writing $a^{\times q}$ instead of $\{a^{\times q}\}$. If multiplicities are $1$, we may drop them from notation, writing $a$ instead of $a^{\times 1}$. 
\end{remark}

\begin{warning}
    For a band $B$ and $q \geq 2$, $B^{\times q}$ refers to $B$ appearing $q$ times in a multi-set while $B^q$ refers to the non-minimal band that is obtained by multiplying $B$ with itself $q$ times as a path (or word). 
\end{warning}

\subsection{Homomorphisms of string and band modules}

There is an equally combinatorial explanation of the homomorphisms between indecomposable string and band modules as shown in \cite{CB_Maps} and \cite{KRAUSE1991}.

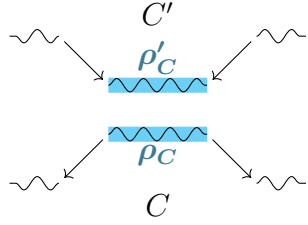
\begin{figure}[H]
    \centering
    \begin{tikzpicture}[scale = 0.65, baseline={([yshift={-\ht\strutbox}]current bounding box.east)}]
    \begin{scope}[name prefix = top-]
        \coordinate (L) at (-2,2){};
        \coordinate (R) at (2,2){};
        \coordinate (LL) at (-3,2){};
        \coordinate (RR) at (3,2){};
    \end{scope}
    \begin{scope}[name prefix = mid-top-]
        \coordinate (L) at (-1,1){};
        \coordinate (R) at (1,1){};
    \end{scope}
    \begin{scope}[name prefix = mid-bot-]
        \coordinate (L) at (-1,0){};
        \coordinate (R) at (1,0){};
    \end{scope}
    \begin{scope}[name prefix = bot-]
        \coordinate (L) at (-2,-1){};
        \coordinate (R) at (2,-1){};
        \coordinate (LL) at (-3,-1){};
        \coordinate (RR) at (3,-1){};
    \end{scope}
    
    \node[anchor = south] at (0,2) {$C'$};
    \node[anchor = north] at (0,-1) {$C$};
    
    \node[anchor = south, color = cyan!50!black] at (0,1) {$\boldsymbol{\rho_C'}$};
    \node[anchor = north, color = cyan!50!black] at (0,0) {$\boldsymbol{\rho_{C}}$};

    \draw[cyan!50!white, line width = 0.2cm, rounded corners] (mid-top-L) -- (mid-top-R);
    \draw[cyan!50!white, line width = 0.2cm, rounded corners] (mid-bot-L) -- (mid-bot-R);

    \draw[->, shorten <=0.1cm, shorten >=0.1cm] (top-L) -- (mid-top-L);
    \draw[->, shorten <=0.1cm, shorten >=0.1cm] (top-R) -- (mid-top-R);
    \draw[->, shorten <=0.1cm, shorten >=0.1cm] (mid-bot-L) -- (bot-L);
    \draw[->, shorten <=0.1cm, shorten >=0.1cm] (mid-bot-R) -- (bot-R);

    \draw[snake it] (mid-top-L) -- (mid-top-R);
    \draw[snake it] (mid-bot-L) -- (mid-bot-R);
    \draw[snake it] (top-L) -- (top-LL);
    \draw[snake it] (top-R) -- (top-RR);
    \draw[snake it] (bot-L) -- (bot-LL);
    \draw[snake it] (bot-R) -- (bot-RR);
\end{tikzpicture}
    \caption{Configuration of substrings giving rise to basis homomorphisms.}
    \label{fig:hom}
\end{figure}

Let $C$ and $C'$ be strings. Let $\rho = (\rho_C, \rho_{C'})$ be a pair in $\Sigma_\text{top}^E(C) \times \Sigma_\text{bot}^E(C')$ for some string $E$, see \Cref{fig:hom}. 
Every such pair gives rise to a homomorphism from $M(C)$ to $M(C')$ by essentially sending vectors coming from $\rho_C$ to their partners in $\rho_{C'}$ and all other elements to zero. A detailed description of these homomorphisms as well as analogues for band modules can be found in \cite{KRAUSE1991}. These morphisms give rise to a $K$-basis of the homomorphism spaces between indecomposable modules. For the purpose of the following, we will only need to consider the size of the $K$-bases of homomorphism spaces from a string module to a string or band module.

\begin{theorem}[\cite{KRAUSE1991}]
\label{thm:hom}
    Let $C$ and $C'$ be strings, $B$ a minimal band, $\lambda \in K^*$ and $q \geq 1$. 
    \begin{itemize}
        \item The homomorphism space $\Hom(M(C),M(C'))$ has a canonical basis in bijection with
        \begin{equation*}
            \bigsqcup_{E \in \strings} \Sigma_\text{top}^E(C) \times \Sigma_\text{bot}^E(C').
        \end{equation*}
        \item The homomorphism space $\Hom(M(C),M(B, \lambda, q))$ has a canonical basis in bijection to 
        \begin{equation*}
            \bigsqcup_{1 \leq i \leq q} \bigsqcup_{E \in \strings} \Sigma_\text{top}^E(C) \times \Sigma_\text{bot}^E(B^\infty).
        \end{equation*}
    \end{itemize}
\end{theorem}

\subsection{Module varieties}

Let $\dd = (d_1, \ldots, d_n) \in \bbN^n$ be a dimension vector. Up to isomorphism, every $\dd$-dimensional module can be considered over the graded $K$-vector space $K^\dd$ that is $M_i = K^{d_i}$ for every $i \in Q_0$. Then every such module can be seen as a $K$-rational point in the affine space 
\begin{equation*}
    \mod(Q,\dd) := \prod_{a \in Q_1} \Hom_K (K^{d_{s(a)}}, K^{d_{t(a)}}).
\end{equation*}
We let $\mod (A,\dd)$ be the affine scheme of representations of $A$ with dimension vector $\dd$ which we regard as the Zariski closed subset of $\mod(Q,\dd)$ of points which satisfy the relations in $I$. From now on, we identify the module with the $K$-rational point.

There is a group action of $\GLdd = \GL_{d_1} \times \dots \times \GL_{d_n}$ on $\mod(A,\dd)$ given by 
\begin{equation*}
    g.M := (M_i, g_{t(a)}M_a g_{s(a)}^{-1})_{i \in Q_0,\, a \in Q_1}
\end{equation*}
where $g = (g_1, \ldots, g_n) \in \GLdd$ and $M \in \mod(A,\dd)$. We denote the orbit of $M$ by $\sO_M$. Note that a module $N$ is in the orbit of $M$ if and only if $N \cong M$.

\biblio

\section{Families of strings and bands}
\label{sect:families}

We divide $\mod(A,\dd)$ into pairwise disjoint constructible subsets which are the unions of module orbits based on the diagram of the modules. In \Cref{sect:degeneration}, we will describe an order of these subsets based on their closures.

\begin{definition} 
\label[definition]{def:families_min_bands}
    Let $\sD$ be a multi-set of strings and minimal bands. We define the \df{family} of $\sD$ as
    \begin{equation*}
        \sO_\sD = \bigcup_{\substack{M \in \mod(A) \\ \diagram{M} = \sD}} \sO_M.
    \end{equation*}
    We define the \df{dimension vector} of $\sD$ as $\ddim(\sD) := \ddim(M)$ for any $M \in \sO_\sD$.
\end{definition}

\begin{remark}
\label[remark]{rmk:constructible}
    Note that any familiy $\sO_\sD$ is a constructible subset of $\mod(A,\dd)$.
\end{remark}

\begin{remark}
\label[remark]{rmk:families}
    All modules in $\sO_\sD$ have the same dimension vector. So the dimension vector of $\sD$ is well-defined. In particular, $\sO_\sD \subseteq \mod(A,\ddim(\sD))$. 
\end{remark}

\begin{remark}
\label[remark]{rmk:families_string_alg}
    Note that this definition can be more broadly applied to string algebras where modules are also given by strings and bands. 
\end{remark}

By definition, these families are disjoint and cover the whole module variety.
\begin{lemma}
    Let $\dd$ be a dimension vector. Then
    \begin{equation*}
        \mod(A,\dd) = \bigsqcup_{\substack{\sD \in \minmsets\\\ddim(\sD) = \dd}} \sO_\sD
    \end{equation*}
    where $\sqcup$ denotes the disjoint union.
\end{lemma}

\begin{remark}
\label[remark]{rmk:partitions}
    For a minimal band $B$ and a multiplicity $q$, the familiy of $B^{\times q}$ can also be described as
    \begin{equation*}
        \sO_{B^{\times q}} := 
        \bigcup_{\substack{q_1 + \ldots + q_\ell = q \\ \text{partition of } q}} \;
        \bigcup_{\lambda_i \in K^*}\sO_{M(B,\lambda_1,q_1) \oplus \ldots \oplus M(B,\lambda_\ell,q_\ell)}.
    \end{equation*}
\end{remark}

\begin{definition}
\label[definition]{def:orbit_non_minimal}
    Let $B^t$ be a non-minimal band and $q$ a multiplicity. We define the \df{family} of $(B^t)^{\times q}$ as
    \begin{equation*}
        \sO_{(B^t)^{\times q}} := 
        \bigcup_{\substack{q_1 + \ldots + q_\ell = q \\ \text{partition of } q}} \;
        \bigcup_{\lambda_i \in K^*}\sO_{M(B^t,\lambda_1,q_1) \oplus \ldots \oplus M(B^t,\lambda_\ell,q_\ell)}.
    \end{equation*}
    Let $\sD = \left\{ D_1^{\times q_1}, \ldots, D_\ell^{\times q_\ell} \right\}$ be a multi-set of (not necessarily minimal) bands and strings where the $D_i$ are pairwise different. We define the \df{family} of $\sD$ as
    \begin{equation*}
        \sO_\sD = \bigcup_{M_i \in \sO_{D_i}} \sO_{M_1 \oplus \ldots \oplus M_\ell}.
    \end{equation*}
    For multi-sets of minimal bands and strings, this definition coincides with \Cref{def:families_min_bands}.
\end{definition}

If a multi-set of strings and bands is of the form $B$ or $C$ for a band or string, we will often talk of the band or string orbit. 

\begin{example}
\label[example]{ex:families}
    Let $A$ be the path algebra of the Kronecker quiver 
$\begin{tikzcd}
   1 & 2.
   \arrow["a", shift left, from=1-1, to=1-2]
   \arrow["b"', shift right, from=1-1, to=1-2]
\end{tikzcd}$ 
We consider the module variety associated to the dimension vector $\dd = (1,1)$. It is isomorphic to $K^2$ equipped with the Zariski topology. The modules of this dimension up to isomorphism are $S_1 \oplus S_2$, $M(a)$, $M(b)$ and $M(ba^-, \lambda, 1)$ for $\lambda \in K^*$. The orbits of these modules are the origin or lines without the origin and can be seen schematically in \Cref{fig:Kronecker_orbits}. The orbit of $S_1 \oplus S_2$ is the origin. The thick lines show the orbit of $M(a)$ and $M(b)$. The thin lines represent the individual orbits of the band modules. We note that closing these module orbits is the same as adding the origin. In contrast, the familiy of the band $B = ba^-$, depicted in \Cref{fig:Kronecker_families} as the coloured area, includes all of $\mod(A,\dd)$ in its closure.
\begin{figure}[H]
    \centering
    \begin{minipage}{.5\textwidth}
      \centering
      \begin{tikzpicture}[scale = 0.75]
            \draw node[] at (0,0){$(0,0)$};
            \draw[line width=1.5pt, shorten <=0.5cm] (0,0) -- (3,0);
            \draw[line width=1.5pt, shorten <=0.5cm] (0,0) -- (0,3);
            \draw[line width = 1.5pt, shorten <=0.5cm] (0,0) -- (0,-3);
            \draw[line width = 1.5pt, shorten <=0.5cm] (0,0) -- (-3,0);
            \foreach \phi in {0,15,...,360}{
                \draw[shorten <=0.5cm] (0,0) -- (\phi:3);
            }
        \end{tikzpicture}
    \caption{Module orbits for the dimension vector $(1,1)$.}
    \label{fig:Kronecker_orbits}
    \end{minipage}%
    \begin{minipage}{.5\textwidth}
      \centering
      \begin{tikzpicture}[scale = 0.75]
            \filldraw[cyan!50!white, opacity = 0.5] (3,3) -- (3,-3) -- (-3,-3) -- (-3,3);
            \draw node[circle,draw=white, fill=white, inner sep=0pt, minimum size=0.5cm] at (0,0){$(0,0)$};
            \draw[line width=1.5pt, shorten <=0.5cm, draw = white] (0,0) -- (3,0);
            \draw[line width=1.5pt, shorten <=0.5cm, draw = white] (0,0) -- (0,3);
            \draw[line width = 1.5pt, shorten <=0.5cm, draw = white] (0,0) -- (0,-3);
            \draw[line width = 1.5pt, shorten <=0.5cm, draw = white] (0,0) -- (-3,0);
        \end{tikzpicture}
    \caption{Family of the band $ba^-$.}
    \label{fig:Kronecker_families}
    \end{minipage}
\end{figure}
\end{example}

The families of \Cref{def:families_min_bands} are also described by a numerical function. In \cite{GLS_gentle_2021}, the authors describe the irreducible components of $\mod(A,\dd)$ via rank functions. A rank function on $\mod(A,\dd)$ is a vector $\rr_M = (\rk(M_a))_{a\in Q_1} \in \bbN^{Q_1}$ for some $M \in \mod(A,\dd)$ where $\rk(M_a)$ is the rank of the linear map $M_a$. For a rank function $\rr$, we define 
\begin{equation*}
    \mod(A,\dd)_{\leq \rr} = \{ M \in \mod(A,\dd) \mid \forall a \in Q_1, \rk(M_a) \leq \rr_a \}.
\end{equation*}
Note that in \cite{GLS_gentle_2021}, this subset is denoted by $\mod(A,\dd,\rr)$. 

\begin{theorem}[\cite{GLS_gentle_2021}] 
\label{thm:irreducible_comp}
    The irreducible components of $\mod(A,\dd)$ are exactly the $\mod(A,\dd)_{\leq \rr}$ where $\rr$ is a maximal rank function with regard to the entry-wise order on $\bbN^{Q_1}$.
\end{theorem}

We will now introduce a refinement of the $\rr_M$ that describe the families defined above. Let 
\begin{equation*}
    h: \mod(A,\dd) \rightarrow \bbN^\strings, M~ \mapsto (\hom(M(C),M))_{C \in \strings}
\end{equation*}
where $\hom(-,?) = \dim_K \Hom(-,?)$ and $\strings$ is the set of all strings of $A$ up to equivalence described in \Cref{sect:preliminaries}. For $\hh \in \bbN^\strings$, we define 
\begin{align*}
    \mod(A,\dd,\hh) &= \{ M \in \mod(A,\dd) \mid h(M) = \hh \},\quad \text{and}\\
    \mod(A,\dd,\geq \hh) &= \{ M \in \mod(A,\dd) \mid \forall C \in \strings, h(M)_C \geq \hh_C \}.
\end{align*}

The following lemma follows directly from \Cref{thm:hom}.

\begin{lemma}
\label[lemma]{lem:h_constant_on_families}
    Let $\sD$ be a multi-set of strings and minimal bands. Then $h$ is constant on $\sO_\sD$. 
\end{lemma}

In particular, it is enough to consider band modules of quasi-length $1$ to understand $h$. By abuse of notation, we may also speak of $h(\sD)$.

First, let us show that $h$ is a refinement of the rank function in the sense that we can derive $\rr_M$ from $h(M)$. Let $a \in Q_1$. We define $\rho_a$ as the maximum length path in $Q$ (not $Q^\pm$) starting at $s(a)$ that does not include $a$.
\begin{equation*}
    \begin{tikzcd}[sep = small]
	&&& {s(a)} \\
	&&&& \bullet \\
	\\
	\bullet
	\arrow["a", from=1-4, to=2-5]
	\arrow["{\rho_a}"', dashed, from=1-4, to=4-1]
\end{tikzcd}
\end{equation*}
For any vertex $i$ in $Q_0$, we define $\rho_i$ as the string 
\begin{equation*}
    \rho_i = \left\{  
    \begin{aligned}
        \rho_a \rho_b^- &\text{ if } \delta^+(i) = \{a,b\}, \\
        \rho_a &\text{ if } \delta^+(i) = \{a\}, \\
        e_i &\text{ if } \delta^+(i) = \varnothing
    \end{aligned}
    \right.
\end{equation*}
where $\delta^+(i)$ is the set of outgoing arrows from $i$.

\begin{lemma}
\label[lemma]{lem:refinement_r_h}
    Let $M \in \mod(A,\dd)$. Then for every $a \in Q_1$,
    \begin{equation*}
        \rk(M_a) = h(M)_{\rho_{s(a)}} - h(M)_{\rho_a}.
    \end{equation*}
\end{lemma}

\begin{proof}
First, let us note that for any $i \in Q_0$, $M(\rho_i) = P_i$ and for any module $M$, $\dim(M_i) = \hom(P_i,M)$.
Without loss of generality, let $M$ be indecomposable. Let us assume that $M \cong M(C)$ for a string $C$.
Note that 
\begin{equation*}
    \hom(M(\rho_a),M) = \sum_{E \in \Sigma_\text{top}(\rho_a)} |\Sigma_\text{bot}^E(C)|.
\end{equation*}
The substrings at the top of $\rho_a$ are exactly the $\rho_{\geq i}$ with $1 \leq i \leq \ell(\rho_a)$. Such a string appears at the bottom of $C$ if and only if there is an appearance of $s(a)$ without the arrow $a$. So
\begin{equation*}
    \hom(M(\rho_a),M) = |\{s(a) \text{ in } C\}| - |\{ a \text{ in } C \}| = \dim M_{s(a)} - \rk_a(M)
\end{equation*} 
This argument directly generalises to the case of band modules.
\end{proof}

\begin{corollary}
\label[corollary]{cor:refinement_r_h}
    Let $M$ be a module. Then $\mod(A,\dd,\geq h(M)) \subseteq \mod(A,\dd)_{\leq \rr_M}$.
\end{corollary}

\begin{proof}
    Follows from \Cref{lem:refinement_r_h} by using that for any $a\in Q_1$, $h(M)_{\rho_{s(a)}} = \dim(M_{s(a)}) = \dd_{s(a)}$.
\end{proof}

The function $h$ has a combinatorial analogue $h'$. 

\begin{definition}
\label[definition]{def:h'}
    Let $\sD = \{C_1, \ldots, C_s, B_1, \ldots, B_t \}$ be a multi-set of strings and minimal bands. We define a function $h': \minmsets \rightarrow \bbN^\strings$ by 
    \begin{equation*}
        h'(\sD) = h'(C_1) + \ldots + h'(C_s) + h'(B_1) + \ldots + h'(B_t)
    \end{equation*}
    where for a band or string $D$, 
    \begin{equation*}
        h'(D) = (|\Sigma_\text{bot}^C(D)|)_{C \in \strings}.
    \end{equation*}
    We also define the function of the same name on the module variety, $h': \mod(A,\dd) \rightarrow \bbN^\strings$ by 
    \begin{equation*}
        h'(M) = h'(\sD)
    \end{equation*}
    if $M \in \sO_\sD$.
\end{definition}

We note that $\hom(M(C),M) = \sum_{E \in \Sigma_\text{top}(C)} h'(M)_E$, so $h$ can be derived from $h'$. The converse is also true.

\begin{lemma}
\label[lemma]{lem:h'_from_h}
    We can inductively derive $h'$ from $h$ by the following formula:
    \begin{equation*}
        h'(\sD)_C = h(\sD)_C - \sum_{E \in \Sigma_\text{top}(C), \ell(E) < \ell(C)} h'(\sD)_E.
    \end{equation*}
\end{lemma}

\begin{proof}
    It is enough to show this for individual string and quasi-length 1 band modules.
    
    Let $M$ be a string module $M(D)$ or a band module $M(D, \lambda, 1)$ for some $\lambda$. 
    Let $C \in \strings$ with $\ell(C) = 0$. Then $C = e_i$ for some $i \in Q_0$. Since $e_i$ is the only substring of $e_i$, $\hom(M(C),M) = h'(D)_C$.
    
    Now let $C \in \strings$ with $\ell(C) > 0$. Then $C$ is its own unique substring on top of itself of length $\ell(C)$. All other substrings have smaller length. So
    \begin{align*}
        \hom(M(C),M) &= \sum_{E \in \Sigma_\text{top}(C)} h'(D)_E \\
        &= h'(D)_C + \sum_{E \in \Sigma_\text{top}(C), \ell(E) < \ell(C)} h'(D)_E.
    \end{align*}
    It follows that 
    \begin{equation*}
        h'(D)_C = \hom(M(C),M) - \sum_{E \in \Sigma_\text{top}(C), \ell(E) < \ell(C)} h'(D)_E.
    \end{equation*}
\end{proof}

\begin{lemma}
\label[lemma]{lem:h'_injective}
    The function $h': \minmsets \rightarrow \bbN^\strings$ is injective.
\end{lemma}

\begin{proof}
    Let $\sD = \{C_1, \ldots, C_s, B_1, \ldots, B_t \}$ and $\sD' = \{C'_1, \ldots, C'_{s'}, B'_1, \ldots, B'_{t'} \}$ be two multi-sets of strings and bands such that $h'(\sD) = h'(\sD')$.
    
    Let $m$ be an integer such that $(B_1)^m$ is longer than any string in $\sD$ and $\sD'$ and such that for any band $B'$ in $\sD$ or $\sD'$ distinct from $B_1$, $m$ satisfies \Cref{lem:bands_substring_length}, i.e. $\forall n \geq m$, for every representative $\wB$ of $B_1$, $(\wB)^n \notin \Sigma_{\text{bot}}((B')^\infty)$. Let $B$ be a representative of $B_1$ such that $B \in \Sigma_{\text{bot}}B_1^\infty$, see \Cref{rmk:band_rota_top_bot}. Then we have
    \begin{equation*}
        \sum_{\substack{B_i \in \sD \\ B_i \sim B_1}} 1 = h'(\sD)_{B^m} = h(\sD')_{B^m} = \sum_{\substack{B'_i \in \sD' \\ B'_i \sim B_1}} 1,
    \end{equation*}
    since no string can have a substring of greater length and no distinct band from $B_1$ can have $B^m$ as a substring at the bottom.
    So $B_1$ appears with the same multiplicity in both multi-sets. Since for two multi-sets $\sT,\sT'$, $h'(\sT \sqcup \sT') = h'(\sT) + h'(\sT')$, we may remove all copies from $B_1$ from the two subsets and continue recursively until all bands are removed from $\sD$ and $\sD'$. 
    
    Let $C_1$ be the string of maximal length with regard to the other strings in $\sD$. It follows that
    \begin{equation*}
        \forall C \in \strings, \ell(C) > \ell(C_1) \Rightarrow h'(\sD)_C = 0.        
    \end{equation*}
    In particular, the length of $C_1$ is also an upper bound for the string-length in $\sD'$. Moreover, for any string $C$, the only substring of $C$ with length $\ell(C)$ is $C$ itself. It follows that 
    \begin{equation*}
        \sum_{\substack{C_i \in \sD \\ C_i \sim C_1}} 1 = h'(\sD)_{C_1} = h(\sD')_{C_1} = \sum_{\substack{C'_i \in \sD' \\ C'_i \sim C_1}} 1.
    \end{equation*}
    So $C_1$ appears with the same multiplicity in both mutli-sets. We may remove all copies of $C_1$ from both multi-sets and continue recursively until both multi-sets are empty and thus the original multi-set $\sD = \sD'$.
\end{proof}

\begin{corollary}
\label[corollary]{cor:families_h}
    For every $\sD \in \minmsets$ of dimension $\dd$, there exists a unique $\hh \in \bbN^\strings$ such that $\sO_\sD = \mod(A,\dd,\hh) $. Moreover, every non-empty $\mod(A,\dd,\hh)$ arises this way.
\end{corollary}

\begin{remark}
    This description also offers a way to see that the $\sO_\sD$ are constructible since $\hom(M, -)$ is upper-semicontinuous, see \cite{CBS2002irreducible}.
\end{remark}

\begin{remark}
\label[remark]{rmk:h_vector_string_alg}
    Note that this result also holds for string algebras.
\end{remark}
    
\biblio{}

\section{Degenerations of strings and bands}
\label{sect:degeneration}

Having defined families in the section before, we may now study their closure. We first note the following.
\begin{lemma}
\label[lemma]{lem:family_closure_quasi_length_1}
    Let $B$ be a band and $q \in \bbN$. Then 
    \begin{equation*}
        \overline{\sO_{B^{\times q}}} = \overline{\bigcup_{\lambda_i \in K^*}\sO_{M(B,\lambda_1,1) \oplus \ldots \oplus M(B,\lambda_q,1)}}.
    \end{equation*}
    In other words, the closure of the family is already obtained as the closure of the union of the direct sum of band modules of quasi-length 1.
\end{lemma}

\begin{proof}
    As per \Cref{rmk:partitions}, the modules in $\sO_{B^{\times q}}$ are of the form \begin{equation*}
        {M(B,\lambda_1,q_1) \oplus \ldots \oplus M(B,\lambda_\ell,q_\ell)}
    \end{equation*} 
    where $q_1 + \ldots + q_\ell = q$ is a partition of $q$ and $\lambda_i \in K^*$.
    These modules only differ at the morphism which contains the $\lambda$ information. In particular, a module of the above form will be determined by a matrix with Jordan blocks $J(\lambda_i,q_i)$. In the variety of square matrices, the family of invertible diagonalizable matrices is a Zariski-dense subset of the family of invertible  matrices. So the closure of the set of matrices with non-zero eigenvalues is already obtained as the closure of the union of the families of invertible diagonalizable matrices. These correspond to direct sums of rank 1 band modules.
\end{proof}

This allows us to show that the families of multi-sets of minimal bands and strings are irreducible.

\begin{lemma}
\label[lemma]{lem:families_irreducible}
    Let $\sD$ be a diagram in $\minmsets$. Then $\sO_\sD$ is irreducible in $\mod(A,\dd)$.
\end{lemma}

\begin{proof}
    Let us first assume that $\sD$ contains no bands, i.e. $\sD = \{ C_1 , \ldots, C_s\}$ where the $C_i$ are strings. Then we can construct a morphism between varieties given by
    \begin{align*}
        \GL_\dd &\mapsto \mod(A,\dd)\\
        g &\mapsto g.(M(C_1) \oplus \ldots \oplus M(C_s)).
    \end{align*}
    This morphism surjects onto $\sO_\sD$. As the variety on the left hand side is irreducible, so is its image $\sO_\sD$.
    Now let us assume that the multi-set is ${\sD = \{ C_1 , \ldots C_s, B_1, \ldots, B_t \}}$ where the $C_i$ are strings and the $B_i$ are minimal bands. We construct a morphism 
    \begin{align*}
        \GL_\dd \times (\bbA \backslash \{0\})^t &\mapsto \mod(A,\dd)\\
        (g, \lambda_1, \ldots, \lambda_t) &\mapsto g.(M(C_1) \oplus \ldots \oplus M(C_s) \oplus M(B_1,\lambda_1,1) \oplus \ldots \oplus M(B_t,\lambda_t,1)).
    \end{align*}
    Note that the variety on the left hand side is a product of irreducible varieties and thus itself irreducible. So its image is also irreducible. This image is the dense subset 
    \begin{equation*}
        \bigcup_{\lambda_i \in K^*}\sO_{M(C_1) \oplus \ldots \oplus M(C_s) \oplus M(B_1,\lambda_1,1) \oplus \ldots \oplus M(B_t,\lambda_t,1)}
    \end{equation*}
    of $\sO_\sD$. Since the closure of an irreducible subset is irreducible, so is $\sO_\sD$.
\end{proof}

\begin{definition}
\label[definition]{def:degeneration_of_strings_and_bands}
    Let $D$ and $D'$ be two multi-sets of strings and bands. Then $\sD$ \textit{degenerates to} $\sD'$, or $\sD'$ is a \textit{degeneration of} $\sD$, denoted $\sD \leqdegfam \sD'$, if $\sO_{\sD'} \subseteq \overline{\sO_\sD}$.
\end{definition}

\begin{lemma}
\label[lemma]{lem:partial_order} 
    Degenerations form a partial order $\leqdegfam$ on the set of multi-sets in $\minmsets$ of dimension $\dd$.
\end{lemma}

\begin{proof}
    Since reflexivity and transitivity follow immediately from the definition, it suffices to show anti-symmetry. We will argue over dimension. Let us take $\sD, \sD' \in \minmsets$. We assume $\sD \leqdeg \sD' \leqdeg \sD$. If we assume that $\sD \neq \sD'$, we obtain that $\sO_{\sD'} \subseteq \overline{\sO_\sD}\backslash \sO_\sD$. Since the families are irreducible, see \Cref{lem:families_irreducible}, this implies $\dim(\sD) < \dim (\sD')$. Similarly, we obtain $\dim(\sD) < \dim (\sD')$. A contradiction.
\end{proof}

This partial order is stable under disjoint union.
\begin{lemma}
\label[lemma]{lem:disjoint_union}
    Let $\sD, \sD', \sE$ be finite multi-sets of strings and bands. Then $\sD \leqdegfam \sD'$ implies $\sD \sqcup \sE \leqdegfam \sD'\sqcup \sE$.
\end{lemma}

We return to the function $h$ described in \Cref{sect:families}. Due to \Cref{cor:families_h} and the fact that the families of diagrams are pairwise disjoint, $h(\sD) = h(\sD')$ if and only of $\sD = \sD'$. This induces a partial order on $\minmsets$ defined by $\sD \leq_h \sD'$ if $h(\sD) \leq h(\sD')$, i.e. if for any string $C$ $h(\sD)_C \leqdeg h(\sD')_C$.

\begin{lemma}
\label[lemma]{lem:two_partial_orders}
    For $\sD, \sD' \in \minmsets$, $\sD \leqdegfam \sD'$ implies that $\sD \leq_h \sD'$.
\end{lemma}

\begin{proof}
    By definition, $\sD \leqdegfam \sD'$ implies that $\sO_{\sD'} \subseteq \overline{\sO_\sD}$. Since $h$ is upper-semicontinuous, it follows that $h(\sD) \leq h(\sD')$. In particular, $\sD \leq_h \sD'$
\end{proof}

\begin{remark}
    It is currently an open question if these partial orders coincide.
\end{remark}

\begin{warning}
    Recall that the classical degeneration order on modules is defined by $M \leqdeg M'$ if $\sO_{M'} \subseteq \overline{\sO_M}$. If we only consider strings, degeneration of their multi-sets will coincide with the degeneration of the associated modules. This is not at all the case if we also consider bands. There is a difference between a degeneration of individual band modules in the classical sense, $M(B, \lambda, 1) \leqdeg M(B', \lambda', 1)$, and a degeneration of bands, $B \leqdeg B'$, in the sense of \Cref{def:degeneration_of_strings_and_bands}. For example, the following \Cref{lem:non_minimal_band_deg} does not hold for degenerations of modules in general. There are also degenerations between strings and bands that do not appear on the module level, see \Cref{ex:band_and_string_deg_vs_module_deg}.
\end{warning}

\begin{lemma}
\label[lemma]{lem:non_minimal_band_deg}
    Let $B$ be a minimal band and $q>1$. Then $B^{\times q} \leqdegfam B^q$.
\end{lemma}

\begin{proof}
    Since $B^q$ is not minimal, the module $M(B^q, \lambda, 1)$ is not indecomposable. We can check that $M(B^q, \lambda, 1) \cong \bigoplus_{\xi \in \bbU_q} M(B, \xi \mu, 1)$ where $\mu$ is a $q$-th root of $\lambda$ and $\bbU_q$ the set of $q$-th roots of unity.
\end{proof}

\begin{example}
    \label[example]{ex:band_and_string_deg_vs_module_deg}
    Let us consider \Cref{ex:families}. To shorten notation, we denote the arrow $b$ by $\drawStringOrBandKronecker{0}{b}{1}$ and $a^-$ by $\drawStringOrBandKronecker{0}{a}{1}$. Bands are marked by circles at the beginning and end, e.g. $\drawStringOrBandKronecker{1}{a,b}{2}$. There is no degeneration between any of the string and band modules with exception of the direct sum of simples, see \Cref{fig:Kronecker_orbits}. However, seen as multi-sets (cf. \cref{fig:Kronecker_families}), the strings $a$ and $b$ both degenerate to the band $\drawStringOrBandKronecker{1}{b,a}{2}$ per \Cref{def:degeneration_of_strings_and_bands}, so $\drawStringOrBandKronecker{1}{b,a}{2} \leqdegfam \drawStringOrBandKronecker{0}{a}{1}$ and $\drawStringOrBandKronecker{1}{b,a}{2} \leqdegfam \drawStringOrBandKronecker{0}{b}{1}$, see \Cref{fig:Kronecker_dim11}. Here, the degeneration order coincides with the $h$-order. While $h$ has infinite length, only finitely many entries are needed to distinguish the order. We write the entries of the strings $(e_1,e_2,a,b)$ in this order underneath the diagrams.
    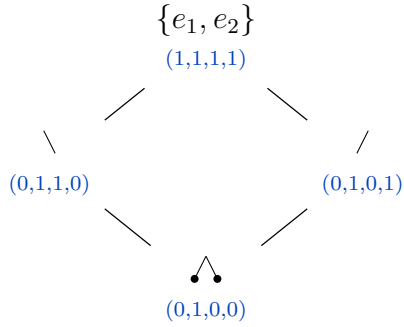
\begin{figure}[H]
        \centering
        \[\begin{tikzcd}[sep = small]
        	& {\stackquiver{\{e_1,e_2\} \\ {\scriptstyle\textcolor{teal!50!blue}{(1,1,1,1)}}}} \\
        	{\stackquiver{\drawStringOrBandKronecker{0}{a}{1} \\ {\scriptstyle\textcolor{teal!50!blue}{(0,1,1,0)}}}} && {\stackquiver{\drawStringOrBandKronecker{0}{b}{1} \\ {\scriptstyle\textcolor{teal!50!blue}{(0,1,0,1)}}}} \\
        	& {\stackquiver{\drawStringOrBandKronecker{1}{b,a}{2} \\ {\scriptstyle\textcolor{teal!50!blue}{(0,1,0,0)}}}}
        	\arrow[no head, from=1-2, to=2-1]
        	\arrow[no head, from=1-2, to=2-3]
        	\arrow[no head, from=2-1, to=3-2]
        	\arrow[no head, from=2-3, to=3-2]
        \end{tikzcd}\]
        \caption{The partial order of degenerations of strings and minimal bands in dimension $(1,1)$. The vectors represent the $h$-vector.}
        \label{fig:Kronecker_dim11}
    \end{figure}
\end{example}

This example conveys a general pattern. We can \textit{delete} arrows in strings or bands to create degenerations.

\begin{definition}
\label[definition]{def:delete_arrow_string} Let $C$ be a string, $m = \ell(C)$ and $1 \leq k \leq m$. Then we define the multi-set of strings arising from $C$ by \df{deleting} the $k$-th arrow as $C_{\hat{k}} = \{ C_{\leq k}, C_{\geq k+1} \}$.
\end{definition}

\begin{example}
\label[example]{ex:deletion_string}
    We reconsider the Kronecker algebra and the string $C = ba^-b$, see \Cref{ex:substring}. Then $C_{\hat{1}} = \{e_2, a^-b\}$, $C_{\hat{2}} = \{b,b\}$, $C_{\hat{3}} = \{ba^-, e_1\}$.
\end{example}

\begin{definition}
\label[definition]{def:delete_arrow_band} Let $B = b_1 \ldots b_m$ be a band and $1 \leq k \leq m$. Then we define the one element multi-set of strings arising from $B$ by \df{deleting} the $k$-th arrow as $B_{\hat{k}} = \{ b_{k+1}\ldots b_m b_1 \ldots b_{k-1} \}$.
\end{definition}

\begin{example}
\label[example]{ex:deletion_band} 
    We consider the Kronecker algebra and the band $B= a^-b$, see \Cref{ex:band}. Then $B_{\hat{1}} = a$ and $B_{\hat{2}} = b$, using the shortened notation for one-element multi-sets and $a^- = a$.
\end{example}

Note that for strings and bands, deleting an arrow preserves the dimension vector of the associated multi-sets.

\begin{lemma}
\label[lemma]{lem:deleting_arrow} 
    Let $D$ be a string or band, $m = \ell(B)$ and $1 \leq k \leq m$. Then $\sO_D \leqdeg \sO_{D_{\hat{k}}}$.
\end{lemma}

\begin{proof}
    We need to show that $\sO_{D_{\hat{k}}} \subset \overline{\sO_D}$. Since $D_{\hat{k}}$ consists only of strings, $\sO_{D_{\hat{k}}}$ is a true $\GLdd$-orbit. So it is enough to show that one module in $\sO_{D_{\hat{k}}}$ is contained in $\overline{\sO_D}$. Let $F$ be the functor used in the construction of the module(s) assigned to $D$. For every $t \in K$,  we define a module $M_t = (K, f_j)_{i,a_j}$ in the source of $F$ with
    \begin{equation*}
        f_j = \left\{
        \begin{aligned}
            t &\quad \textnormal{ if } j = k, \\
           1 &\quad \textnormal{ else.}
        \end{aligned}
        \right.
    \end{equation*}
    Let $\dd = \ddim(D)$. We consider $K$ as a variety with Zariski topology and define a regular map $\varphi: K \rightarrow \mod(A, \dd), t \mapsto F(M_t)$. Note that $F(M_0) \in \sO_{D_{\hat{k}}}$. For every $t \in K^*$, $F(M_t)$ lies in the family of $D$. This shows that indeed $\sO_{D_{\hat{k}}} \subseteq \overline{\sO_D}$.
\end{proof}

Note that we can always successively delete arrows until we reach the multi-set consisting only of $e_i$s.

\begin{remark}
\label[remark]{rmk:deletion_string_alg}
    This result also holds for string algebras.
\end{remark}

\subsection{Reachings}

We wish to study further links between combinatorial traits of strings and bands and degenerations. To this end, we study degenerations which provide a combinatorial mean of constructing degenerations. The following notion of reachings was introduced in the context of gentle algebras in \cite{PPP_kissing}.

\begin{definition}
\label[definition]{def:reaching}
    Let $C$ and $C'$ be strings. We say that $C$ \df{reaches for} $C'$ if there exist substrings $C[i,j]$ on top of $C$ and $C'[i',j']$ at the bottom of $C'$ such that they are equal as strings and satisfy the swinging arms condition:
    \begin{itemize}
        \item If $i'=1$ then $i \neq 1$.
        \item If $j' = \ell(C')+1$ then $j \neq \ell(C)+1$.
    \end{itemize}
    We call the data $C[i,j] = C'[i',j']$ a \df{reaching}.
    We extend this to reachings between two bands or a band and string by replacing substrings of the string with those of the unraveled band. They will automatically fulfill the swinging arms condition.
\end{definition}

\begin{remark}
The reason it is called swinging arms condition is the following. We call the substring $C_{\leq i}$, resp. $C_{\geq j}$, a left, resp. right arm of $C$ if the substring does not have length 0. We do the same for $C'$. The swinging arms condition now translates to: If $C'$ does not have a left, resp. right arm, then $C$ must have one. So the arm ``swings down'', see \Cref{fig:swinging_arms}. Note that the condition does allow up to two arms on each side, so the image will not always resemble a human with swinging arms.
\begin{figure}[H]
    \centering
    \begin{tikzpicture}[scale = 0.65, baseline={([yshift={-\ht\strutbox}]current bounding box.east)}]
    \begin{scope}[name prefix = top-]
        \node at (0,2) {$C'$};
        \coordinate (L) at (-2,2){};
        \coordinate (R) at (2,2){};
        \coordinate (LL) at (-3,2){};
        \coordinate (RR) at (3,2){};
    \end{scope}
    \begin{scope}[name prefix = mid-top-]
        \coordinate (L) at (-1,1){};
        \coordinate (R) at (1,1){};
    \end{scope}
    \begin{scope}[name prefix = mid-bot-]
        \coordinate (L) at (-1,0){};
        \coordinate (R) at (1,0){};
    \end{scope}
    \begin{scope}[name prefix = bot-]
        \node at (0,-1) {$C$};
        \coordinate (L) at (-2,-1){};
        \coordinate (R) at (2,-1){};
        \coordinate (LL) at (-3,-1){};
        \coordinate (RR) at (3,-1){};
    \end{scope}
    
    \draw[yellow, line width = 0.2cm, rounded corners] (top-LL) -- (top-L) -- (mid-top-L);
    
    \draw[cyan!50!white, line width = 0.2cm, rounded corners] (mid-bot-R) -- (bot-R) -- (bot-RR);
        
    \draw[->, shorten <=0.1cm, shorten >=0.1cm] (top-L) -- (mid-top-L);
    \draw[->, shorten <=0.1cm, shorten >=0.1cm] (mid-bot-R) -- (bot-R);

    \draw[snake it] (mid-top-L) -- (mid-top-R);
    \draw[snake it] (mid-bot-L) -- (mid-bot-R);
    \draw[snake it] (top-L) -- (top-LL);
    \draw[snake it] (bot-R) -- (bot-RR);
\end{tikzpicture}
    \caption{Swinging arms condition.}
    \label{fig:swinging_arms}
\end{figure}
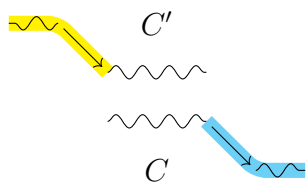
\end{remark}

\begin{example}
\label[example]{ex:reaching}
    We reconsider the Kronecker quiver, see \Cref{ex:band}, and the strings $C = ba^-$ and $C' = e_1$ . Then there is a reaching $C[2,2] = C'[1,1]$.
\end{example}

\begin{remark}
\label[remark]{rmk:reaching}
    Note that $D$ might reach for $D'$ in multiple places, that $D'$ can also reach for $D$ and that $D$ can reach for itself. If $D$ reaches for itself, we formally distinguish between a reaching where we consider $D$ separate from the other copy of $D$, and an \df{auto-reaching} where we consider them as the same object.
\end{remark}

\begin{example}
\label[example]{ex:auto_reaching}
    We reconsider the Kronecker quiver and \Cref{ex:top_string,,ex:bot_string}. The string $C= ba^-b$ reaches for itself with $C_{\leq 2} = b = C_{\geq 3}$.
\end{example}

\subsection{Auto reachings} 

We describe some properties of auto-reachings. 

\begin{lemma}
\label[lemma]{lem:length_band_autoreaching}
    Let $B$ be a band with an auto-reaching given by substrings of the form $L$ of $B^\infty$. Then $\ell(L) \leq \ell(B) -2$.
\end{lemma}

\begin{proof}
    The string $L$ appears as a substring on top  surrounded by $a$ and $b^-$, and as a substring at the bottom surrounded by $c^-$ and $d$ with $a,b,c,d \in Q_1$. 
    Let us assume $\ell(L) \geq \ell(B)$. If the two substrings are oriented in the same direction, the $\ell(B)$-th arrow in $L$ has to be $a$ but also $c^-$ which is impossible. If instead they are oriented in opposite directions, the $\ell(B)$-th arrow in $L$ has to be $a$ but also $d^-$ which remains impossible.
    If $\ell(L) = \ell(B) -1$, $a = b^-$ which also causes a contradiction.
\end{proof}

Due to this length constraint, we can always consider auto-reachings of bands as auto-reachings of string representatives of the band. So the following section on strings can also be applied to them.

We note that there are two cases of auto-reachings: intersecting and non-intersecting ones. 

\begin{definition}
\label[definition]{def:intersecting_auto_reaching_string}
    Let $C$ be a string with an auto-reaching $L = (C[i,j],C[i',j'])$. We say that $L$ is \df{intersecting} if the two substrings overlap in $C$, i.e.
    \begin{equation*}
        i < i' \leq j < j' \qquad \text{or} \qquad i' < i \leq j' < j.
    \end{equation*}
    Note that in particular, $L$ is intersecting if the two substrings intersect in only one vertex.
\end{definition}

\begin{example}
\label[example]{ex:intersecting_reaching_string}
    We consider the Kronecker quiver, see \Cref{ex:auto_reaching}, and the string $C$ represented by $ba^-ba^-b$. Then $C_{\leq 4} = C_{\geq 3}$ defines an intersecting auto-reaching on $C$.
\end{example}

Auto-reachings that are not intersecting are called \df{non-intersecting}. An example was given in \Cref{ex:auto_reaching}. Intersecting auto-reachings have some interesting combinatorial constraints.

\begin{lemma}
\label[lemma]{lem:intersecting_auto_reaching_directionality}
    Let $C$ be a string with an intersecting auto-reaching. Then the two substrings defining the reaching are oriented in the same direction in $C$.
\end{lemma}

\begin{proof}
    Let the reaching be given by the substrings $C[i,j]$ and $C[i',j']$. Since they intersect, they must be of length at least 1. Let us assume that these substrings are oriented in opposite directions, that is the representative $C[i,j]$ is the inverse of the representative $C[i',j']$. If $j = i'$, then the last arrow in $C[i,j]$ would be followed by its own inverse in $C[i',j']$. So let us assume $i' < j$. We consider $C[i',j]$. It is a self-inverse string of length at least 1. But this is impossible as this implies that an arrow in the middle of the string is either immediately followed by its own inverse or is both a normal and an inverse arrow at the same time.
\end{proof}

\begin{lemma}
\label[lemma]{lem:intersecting_auto_reaching_periodicity}
    Let $C$ be a string such that there is an intersecting auto-reaching given by $C[i,j]$ and $C[i',j']$ where $i < i' \leq j < j'$. Then $C[i,i']$ is cyclically equivalent to $C[j,j']$ and $C[i,j]$ is a substring of $(C[i,i'])^\infty$.
\end{lemma}

\begin{proof}
    Let $\ell = \ell(C[i,j])$ and $k = \ell(C[i,i'])$. Let $ \ell = m \cdot k + r$ with $ 0 \leq r \leq k-1$. Now we may use $C[i,j] = C[i',j']$ and the fact that the two substrings intersect on $C[i',j]$ to show that the substring is periodic, i.e. $C[i,j] = (C[i,i'])^m C[i,i+r]$, see \Cref{fig:periodicity_in_auto_reaching}.
    \begin{figure}[H]
        \centering
        \begin{tikzpicture}[scale = 0.5]

    \draw[|-|] (0,2) -- (4,2);  
    \draw[|-|] (4,2) -- (8,2);
    \draw[|-|] (8,2) -- (12,2);
    \draw[|-|] (12,2) -- (15,2);
    \draw[|-|] (4,0) -- (8,0);
    \draw[|-|] (8,0) -- (15,0);
    \draw[|-|] (15,0) -- (19,0);
        
    \fill[pattern color = cyan!50!black, pattern = horizontal lines] (0,2) -- (4,2) -- (8,0) -- (4,0) -- cycle;
    \fill[pattern color = yellow!80!black, pattern = vertical lines] (4,0) -- (4,2) -- (8,2) -- (8,0) -- cycle;
    
    \node[anchor = east] at (0,2) {$C[i,j]$};
    \node[anchor = west] at (19,0) {$C[i',j']$};
    \node[anchor = south] at (2,2) {$C[i,i']$};
    \node[anchor = south, color = yellow!50!black] at (6,2) {$\boldsymbol{C[i,i']}$};
    \node[anchor = north, color = cyan!50!black] at (6,0) {$\boldsymbol{C[i,i']}$};
    \node[anchor = south] at (13.5,2) {$C[i,i+r]$};
    \node[] at (10,1) {$\dots$};
    
\end{tikzpicture}
        \caption{Periodicity in intersecting auto-reachings.}
        \label{fig:periodicity_in_auto_reaching}
    \end{figure}
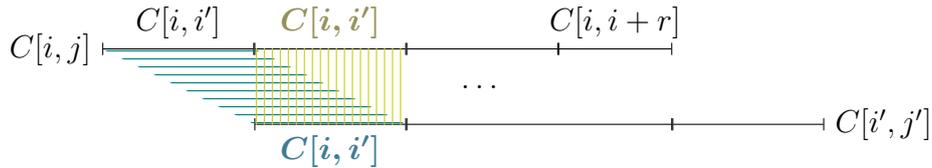
    Since the last $k$ arrows of $C[i,j]$ also define $C[j,j']$, it follows that
    \begin{equation*}
        C[j,j'] = C[i'-(k-r),i'] C[i,i+r].
    \end{equation*} 
    In particular, this is a cyclic permutation of $C[i,i']$.
\end{proof}

\begin{remark}
    Note that in the case of bands, intersecting auto reachings can only intersect in one place not at both ends. Otherwise there would have to be an arrow going in two directions at once, similarly to the proof of \Cref{lem:length_band_autoreaching}.
\end{remark}

\subsection{Resolutions of reachings}

\begin{definition}
\label[definition]{def:reaching_resolution_strings}
    Let $C$, $C'$ be two strings with a reaching $L = C[i,j] = C'[i',j']$. We define the \df{resolution} of the reaching, denoted $\sR_L$, as the multi-set of strings $C_{\leq j}C'_{\geq j'}$ and $C'_{\leq j'}C_{\geq j}$, see \Cref{fig:resolution_strings}.
    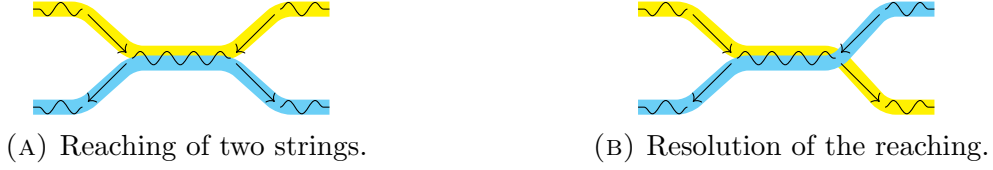
\begin{figure}[H]
    \begin{subfigure}{.5\textwidth}
        \centering
        \begin{tikzpicture}[scale = 0.65, baseline={([yshift={-\ht\strutbox}]current bounding box.east)}]
    \begin{scope}[name prefix = top-]
        \coordinate (L) at (-2,1){};
        \coordinate (R) at (2,1){};
        \coordinate (LL) at (-3,1){};
        \coordinate (RR) at (3,1){};
    \end{scope}
    \begin{scope}[name prefix = mid-]
        \coordinate (L) at (-1,0){};
        \coordinate (R) at (1,0){};
    \end{scope}
    \begin{scope}[name prefix = bot-]
        \coordinate (L) at (-2,-1){};
        \coordinate (R) at (2,-1){};
        \coordinate (LL) at (-3,-1){};
        \coordinate (RR) at (3,-1){};
    \end{scope}

    \draw[yellow, line width = 0.2cm, rounded corners] (top-LL) -- (top-L) -- ($(mid-L)+(0,0.1)$) -- ($(mid-R)+(0,0.1)$) -- (top-R) -- (top-RR);

    \draw[cyan!50!white, line width = 0.2cm, rounded corners] (bot-LL) -- (bot-L) -- ($(mid-L)+(0,-0.1)$) -- ($(mid-R)+(0,-0.1)$) -- (bot-R) -- (bot-RR);
        
    \draw[->, shorten <=0.1cm, shorten >=0.1cm] (top-L) -- (mid-L);
    \draw[->, shorten <=0.1cm, shorten >=0.1cm] (top-R) -- (mid-R);
    \draw[->, shorten <=0.1cm, shorten >=0.1cm] (mid-L) -- (bot-L);
    \draw[->, shorten <=0.1cm, shorten >=0.1cm] (mid-R) -- (bot-R);

    \draw[snake it] (mid-L) -- (mid-R);
    \draw[snake it] (top-L) -- (top-LL);
    \draw[snake it] (top-R) -- (top-RR);
    \draw[snake it] (bot-L) -- (bot-LL);
    \draw[snake it] (bot-R) -- (bot-RR);
\end{tikzpicture}
        \caption{Reaching of two strings.}
    \end{subfigure}%
    \begin{subfigure}{.5\textwidth}
        \centering
        \begin{tikzpicture}[scale = 0.65, baseline={([yshift={-\ht\strutbox}]current bounding box.east)}]
    \begin{scope}[name prefix = top-]
        \coordinate (L) at (-2,1){};
        \coordinate (R) at (2,1){};
        \coordinate (LL) at (-3,1){};
        \coordinate (RR) at (3,1){};
    \end{scope}
    \begin{scope}[name prefix = mid-]
        \coordinate (L) at (-1,0){};
        \coordinate (R) at (1,0){};
    \end{scope}
    \begin{scope}[name prefix = bot-]
        \coordinate (L) at (-2,-1){};
        \coordinate (R) at (2,-1){};
        \coordinate (LL) at (-3,-1){};
        \coordinate (RR) at (3,-1){};
    \end{scope}

    \draw[yellow, line width = 0.2cm, rounded corners] (top-LL) -- (top-L) -- ($(mid-L)+(0,0.1)$) -- ($(mid-R)+(0,0.1)$) -- (bot-R) -- (bot-RR);

    \draw[cyan!50!white, line width = 0.2cm, rounded corners] (bot-LL) -- (bot-L) -- ($(mid-L)+(0,-0.1)$) -- ($(mid-R)+(0,-0.1)$) -- (top-R) -- (top-RR);
        
    \draw[->, shorten <=0.1cm, shorten >=0.1cm] (top-L) -- (mid-L);
    \draw[->, shorten <=0.1cm, shorten >=0.1cm] (top-R) -- (mid-R);
    \draw[->, shorten <=0.1cm, shorten >=0.1cm] (mid-L) -- (bot-L);
    \draw[->, shorten <=0.1cm, shorten >=0.1cm] (mid-R) -- (bot-R);

    \draw[snake it] (mid-L) -- (mid-R);
    \draw[snake it] (top-L) -- (top-LL);
    \draw[snake it] (top-R) -- (top-RR);
    \draw[snake it] (bot-L) -- (bot-LL);
    \draw[snake it] (bot-R) -- (bot-RR);
\end{tikzpicture}
        \caption{Resolution of the reaching.}
    \end{subfigure}
    \caption{Resolution of a reaching between two strings.}
    \label{fig:resolution_strings}
    \end{figure}
\end{definition}

\begin{example}
\label[example]{ex:resolution_reaching_strings}
    We reconsider \Cref{ex:reaching}. The resolution of this reaching is $\{a,b\}$.
\end{example}

\begin{remark}
    Note that the above is always well-defined if the substring defining the reaching has at least length 1. For a reaching in only a vertex, for a chosen representative of $C$, there is exactly one representative of $C'$ such that the above is well-defined and the relations of $A$ are avoided. We always implicitly choose this orientation.
\end{remark}

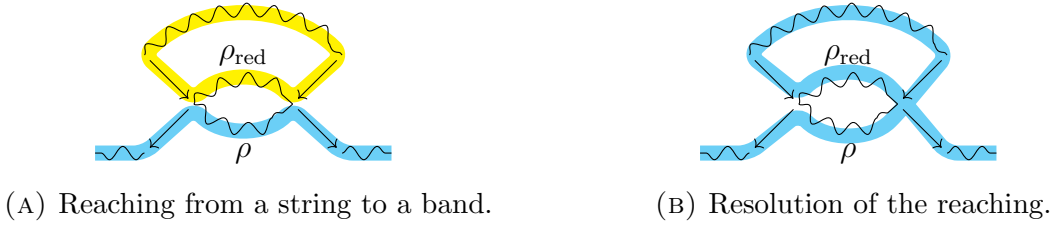
\begin{figure}[H]
    \begin{subfigure}{.5\textwidth}
        \centering
        \begin{tikzpicture}[scale = 0.65, baseline={([yshift={-\ht\strutbox}]current bounding box.east)}]
    \begin{scope}[name prefix = top-]
        \coordinate (L) at (-2,1){};
        \coordinate (R) at (2,1){};
    \end{scope}
    \begin{scope}[name prefix = mid-]
        \coordinate (L) at (-1,0){};
        \coordinate (R) at (1,0){};
    \end{scope}
    \begin{scope}[name prefix = bot-]
        \coordinate (L) at (-2,-1){};
        \coordinate (R) at (2,-1){};
        \coordinate (LL) at (-3,-1){};
        \coordinate (RR) at (3,-1){};
    \end{scope}

    \draw[yellow, line width = 0.2cm, rounded corners]  ($(mid-L)+(0,0.1)$) to[out=45,in=135] ($(mid-R)+(0,0.1)$) -- (top-R) to[out=135,in=45] (top-L) -- cycle;

    \draw[cyan!50!white, line width = 0.2cm, rounded corners] (bot-LL) -- (bot-L) -- ($(mid-L)+(0,-0.1)$) to[out=315,in=225] ($(mid-R)+(0,-0.1)$) -- (bot-R) -- (bot-RR);
        
    \draw[->, shorten <=0.1cm, shorten >=0.1cm] (top-L) -- (mid-L);
    \draw[->, shorten <=0.1cm, shorten >=0.1cm] (top-R) -- (mid-R);
    \draw[->, shorten <=0.1cm, shorten >=0.1cm] (mid-L) -- (bot-L);
    \draw[->, shorten <=0.1cm, shorten >=0.1cm] (mid-R) -- (bot-R);

    \draw[snake it] (mid-L) to[out=45,in=135] (mid-R);
    \draw[snake it] (mid-L) to[out=315,in=225] (mid-R);
    \draw[snake it] (top-L) to[out=45,in=135] (top-R);
    \draw[snake it] (bot-L) -- (bot-LL);
    \draw[snake it] (bot-R) -- (bot-RR);
    
    \node at (0,1) {$\rho_\text{red}$};
    \node at (0,-1) {$\rho$};
    
\end{tikzpicture}
        \caption{Reaching from a string to a band.}
    \end{subfigure}%
    \begin{subfigure}{.5\textwidth}
        \centering
        \begin{tikzpicture}[scale = 0.65, baseline={([yshift={-\ht\strutbox}]current bounding box.east)}]
    \begin{scope}[name prefix = top-]
        \coordinate (L) at (-2,1){};
        \coordinate (R) at (2,1){};
    \end{scope}
    \begin{scope}[name prefix = mid-]
        \coordinate (L) at (-1,0){};
        \coordinate (R) at (1,0){};
    \end{scope}
    \begin{scope}[name prefix = bot-]
        \coordinate (L) at (-2,-1){};
        \coordinate (R) at (2,-1){};
        \coordinate (LL) at (-3,-1){};
        \coordinate (RR) at (3,-1){};
    \end{scope}

    \draw[cyan!50!white, line width = 0.2cm, rounded corners] (bot-LL) -- (bot-L) -- ($(mid-L)+(0,-0.2)$) to[out=315,in=225] ($(mid-R)+(0,-0.2)$) -- (top-R) to[out=135,in=45] (top-L) -- ($(mid-L)+(0,0.2)$) to[out=45,in=135] ($(mid-R)+(0,0.2)$) -- (bot-R) -- (bot-RR);
        
    \draw[->, shorten <=0.1cm, shorten >=0.1cm] (top-L) -- (mid-L);
    \draw[->, shorten <=0.1cm, shorten >=0.1cm] (top-R) -- (mid-R);
    \draw[->, shorten <=0.1cm, shorten >=0.1cm] (mid-L) -- (bot-L);
    \draw[->, shorten <=0.1cm, shorten >=0.1cm] (mid-R) -- (bot-R);

    \draw[snake it] (mid-L) to[out=45,in=135] (mid-R);
    \draw[snake it] (mid-L) to[out=315,in=225] (mid-R);
    \draw[snake it] (top-L) to[out=45,in=135] (top-R);
    \draw[snake it] (bot-L) -- (bot-LL);
    \draw[snake it] (bot-R) -- (bot-RR);
    
    \node at (0,1) {$\rho_\text{red}$};
    \node at (0,-1) {$\rho$};
\end{tikzpicture}
        \caption{Resolution of the reaching.}
    \end{subfigure}
    \caption{Resolution of a reaching from a string to a band.}
    \label{fig:resolution_string_band}
\end{figure}

We wish to extend this definition to bands as well. Let $B$ be a band and $C$ a string such that there is a reaching $L$ between them. So there is a substring $\rho$ at the bottom or on top of $B^\infty$ defining one half of the reaching. Let $\rho_{\text{red}}$ be the substring of $B$ associated to $\rho$, see \Cref{def:substring_band}. Without loss of generality, we choose a representative of $B$ such that $\rho_\text{red} = B[i,j]$ with $2\leq i \leq j \leq \ell(B)$. Then we use the above definition to define a resolution of this reaching, replacing  $\rho$ with $\rho_\text{red}$. Instead of obtaining two strings, we will obtain one string, as we still consider $B$ as a band and respect its structure in the resolution, see \Cref{fig:resolution_string_band}.

\begin{example}
\label[example]{ex:reaching_resolution_string_band}
    We consider the gentle algebra given by quiver and relations 
    \begin{equation*}
       \begin{tikzcd}
        	1 & 2 & 3.
        	\arrow[""{name=0, anchor=center, inner sep=0}, "a", shift left, from=1-1, to=1-2]
        	\arrow[""{name=1, anchor=center, inner sep=0}, "b"', shift right, from=1-1, to=1-2]
        	\arrow[""{name=2, anchor=center, inner sep=0}, "c", shift left, from=1-2, to=1-3]
        	\arrow[""{name=3, anchor=center, inner sep=0}, "d"', shift right, from=1-2, to=1-3]
        	\arrow[curve={height=-12pt}, shorten <=8pt, shorten >=8pt, dashed, no head, from=0, to=2]
        	\arrow[curve={height=12pt}, shorten <=8pt, shorten >=8pt, dashed, no head, from=1, to=3]
        \end{tikzcd}
    \end{equation*}
    For shorter notation, we write
    \begin{equation*}
        ba^- = \drawStringOrBand{0}{b,a}{2}, \qquad d^-c = \drawStringOrBand{0}{d,c}{2}.
    \end{equation*}
    Bands are marked by circles at the beginning and end, e.g. $B = \drawStringOrBand{1}{d,c}{2}$. We consider the string $C = \drawStringOrBand{0}{b,a,d,c,d,c}{6}$. There is a reaching from $B$ to $C$ given by the substring $\rho = \drawStringOrBand{0}{d,c,d,c}{4}$. Then $\rho_\text{red}$ is $e_2$ and the chosen representative for $B$ has to be $\drawStringOrBand{1}{c,d}{2}$. The resolution of the reaching is the string $\drawStringOrBand{0}{b,a,d,c,d,c,d,c}{8}$.
\end{example}    

\begin{figure}[H]
    \begin{subfigure}{.5\textwidth}
        \centering
        \begin{tikzpicture}[scale = 0.65, baseline={([yshift={-\ht\strutbox}]current bounding box.east)}]
    \begin{scope}[name prefix = top-]
        \coordinate (L) at (-2,1){};
        \coordinate (R) at (2,1){};
    \end{scope}
    \begin{scope}[name prefix = mid-]
        \coordinate (L) at (-1,0){};
        \coordinate (R) at (1,0){};
    \end{scope}
    \begin{scope}[name prefix = bot-]
        \coordinate (L) at (-2,-1){};
        \coordinate (R) at (2,-1){};
    \end{scope}

    \draw[yellow, line width = 0.2cm, rounded corners]  ($(mid-L)+(0,0.1)$) to[out=45,in=135] ($(mid-R)+(0,0.1)$) -- (top-R) to[out=135,in=45] (top-L) -- cycle;

    \draw[cyan!50!white, line width = 0.2cm, rounded corners] ($(mid-L)+(0,-0.1)$) to[out=315,in=225] ($(mid-R)+(0,-0.1)$) -- (bot-R) to[out=225,in=315] (bot-L) -- cycle;
        
    \draw[->, shorten <=0.1cm, shorten >=0.1cm] (top-L) -- (mid-L);
    \draw[->, shorten <=0.1cm, shorten >=0.1cm] (top-R) -- (mid-R);
    \draw[->, shorten <=0.1cm, shorten >=0.1cm] (mid-L) -- (bot-L);
    \draw[->, shorten <=0.1cm, shorten >=0.1cm] (mid-R) -- (bot-R);

    \draw[snake it] (mid-L) to[out=45,in=135] (mid-R);
    \draw[snake it] (mid-L) to[out=315,in=225] (mid-R);
    \draw[snake it] (top-L) to[out=45,in=135] (top-R);
    \draw[snake it] (bot-L) to[out=315,in=225] (bot-R);
        
    \node at (0,1) {$\rho'_\text{red}$};
    \node at (0,-1) {$\rho_\text{red}$};
\end{tikzpicture}
        \caption{Reaching of two bands.}
    \end{subfigure}%
    \begin{subfigure}{.5\textwidth}
        \centering
        \begin{tikzpicture}[scale = 0.65, baseline={([yshift={-\ht\strutbox}]current bounding box.east)}]
    \begin{scope}[name prefix = top-]
        \coordinate (L) at (-2,1){};
        \coordinate (R) at (2,1){};
    \end{scope}
    \begin{scope}[name prefix = mid-]
        \coordinate (L) at (-1,0){};
        \coordinate (R) at (1,0){};
    \end{scope}
    \begin{scope}[name prefix = bot-]
        \coordinate (L) at (-2,-1){};
        \coordinate (R) at (2,-1){};
    \end{scope}

    \draw[cyan!50!white, line width = 0.2cm, rounded corners]  ($(mid-L)+(0,0.2)$) to[out=45,in=135] ($(mid-R)+(0,0.2)$) -- (bot-R) to[out=225,in=315] (bot-L) -- ($(mid-L)+(0,-0.2)$) to[out=315,in=225] ($(mid-R)+(0,-0.2)$) -- (top-R) to[out=135,in=45] (top-L) -- cycle;
        
    \draw[->, shorten <=0.1cm, shorten >=0.1cm] (top-L) -- (mid-L);
    \draw[->, shorten <=0.1cm, shorten >=0.1cm] (top-R) -- (mid-R);
    \draw[->, shorten <=0.1cm, shorten >=0.1cm] (mid-L) -- (bot-L);
    \draw[->, shorten <=0.1cm, shorten >=0.1cm] (mid-R) -- (bot-R);

    \draw[snake it] (mid-L) to[out=45,in=135] (mid-R);
    \draw[snake it] (mid-L) to[out=315,in=225] (mid-R);
    \draw[snake it] (top-L) to[out=45,in=135] (top-R);
    \draw[snake it] (bot-L) to[out=315,in=225] (bot-R);
        
    \node at (0,1) {$\rho'_\text{red}$};
    \node at (0,-1) {$\rho_\text{red}$};
\end{tikzpicture}
        \caption{Resolution of the reaching.}
    \end{subfigure}
    \caption{Resolution of a reaching between two bands.}
    \label{fig:resolution_bands}
\end{figure}
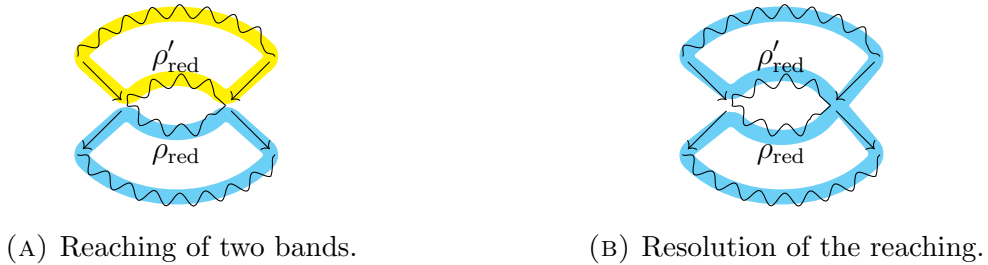

We can extend this definition the same way to encompass reachings of two bands. The resolution will be a single band. Note that this involves up to two reduced substrings: $\rho_\text{red}$ and $\rho'_\text{red}$, see \cref{fig:resolution_bands}.

\begin{remark}
    Note that there is a priori a choice for when $\rho_{red}$ or $\rho'_{red}$ is walked. However, due to the periodicity of $\rho$ as a substring of $B$, the resolution is indeed independent of this choice. It is also independent of the choice of representative of $B$. However, resulting bands are not necessarily minimal even if we start with minimal bands.
\end{remark}

\begin{example}
\label[example]{ex:resolution_produces_non_minimal_bands}
    We consider the algebra from \Cref{ex:reaching_resolution_string_band}. Let $B = \drawStringOrBand{1}{b,a,d,c,d,c}{6}$ and $B' = \drawStringOrBand{1}{a,b}{2}$. Then $B$ reaches for $B'$ via $e_2$. The resolution thereof is $\drawStringOrBand{1}{b,a,d,c,b,a,d,c}{8}$ which is not minimal.
\end{example}

Next, we resolve auto-reachings. We start with the non-intersecting ones. Recall that due to \Cref{lem:length_band_autoreaching}, we can always depict auto-reachings of bands as though they happened on strings.

\begin{definition}
\label[definition]{def:non_intersecting_auto_reaching_resolution}
    Let $D$ be a string or band with a reaching $\rho$. We define the \df{resolution} of $\rho$ as the multi-set $\sR^a_\rho$ similarly to \Cref{def:reaching_resolution_strings}, see \Cref{fig:resolution_auto}, where some of the $\ast_i$ may be identified depending on the directionality of the substrings defining $\rho$ and on $D$ being a band or string. They remain identified after the resolution.
    \begin{figure}[H]
        \begin{subfigure}{.5\textwidth}
            \centering
            \begin{tikzpicture}[scale = 0.65, baseline={([yshift={-\ht\strutbox}]current bounding box.east)}]
    \begin{scope}[name prefix = top-]
        \node[anchor = east] at (-3,1) {$\ast_1$};
        \node[anchor = west] at (3,1) {$\ast_2$};
        \coordinate (L) at (-2,1){};
        \coordinate (R) at (2,1){};
        \coordinate (LL) at (-3,1){};
        \coordinate (RR) at (3,1){};
    \end{scope}
    \begin{scope}[name prefix = mid-]
        \coordinate (L) at (-1,0){};
        \coordinate (R) at (1,0){};
    \end{scope}
    \begin{scope}[name prefix = bot-]
        \node[anchor = east] at (-3,-1) {$\ast_3$};
        \node[anchor = west] at (3,-1) {$\ast_4$};
        \coordinate (L) at (-2,-1){};
        \coordinate (R) at (2,-1){};
        \coordinate (LL) at (-3,-1){};
        \coordinate (RR) at (3,-1){};
    \end{scope}

    \draw[yellow, line width = 0.2cm, rounded corners] (top-LL) -- (top-L) -- ($(mid-L)+(0,0.1)$) -- ($(mid-R)+(0,0.1)$) -- (top-R) -- (top-RR);

    \draw[cyan!50!white, line width = 0.2cm, rounded corners] (bot-LL) -- (bot-L) -- ($(mid-L)+(0,-0.1)$) -- ($(mid-R)+(0,-0.1)$) -- (bot-R) -- (bot-RR);
        
    \draw[->, shorten <=0.1cm, shorten >=0.1cm] (top-L) -- (mid-L);
    \draw[->, shorten <=0.1cm, shorten >=0.1cm] (top-R) -- (mid-R);
    \draw[->, shorten <=0.1cm, shorten >=0.1cm] (mid-L) -- (bot-L);
    \draw[->, shorten <=0.1cm, shorten >=0.1cm] (mid-R) -- (bot-R);

    \draw[snake it] (mid-L) -- (mid-R);
    \draw[snake it] (top-L) -- (top-LL);
    \draw[snake it] (top-R) -- (top-RR);
    \draw[snake it] (bot-L) -- (bot-LL);
    \draw[snake it] (bot-R) -- (bot-RR);
\end{tikzpicture}
            \caption{Non-intersecting auto-reaching.}
        \end{subfigure}%
        \begin{subfigure}{.5\textwidth}
            \centering
            \begin{tikzpicture}[scale = 0.65, baseline={([yshift={-\ht\strutbox}]current bounding box.east)}]
    \begin{scope}[name prefix = top-]
        \node[anchor = east] at (-3,1) {$\ast_1$};
        \node[anchor = west] at (3,1) {$\ast_2$};
        \coordinate (L) at (-2,1){};
        \coordinate (R) at (2,1){};
        \coordinate (LL) at (-3,1){};
        \coordinate (RR) at (3,1){};
    \end{scope}
    \begin{scope}[name prefix = mid-]
        \coordinate (L) at (-1,0){};
        \coordinate (R) at (1,0){};
    \end{scope}
    \begin{scope}[name prefix = bot-]
        \node[anchor = east] at (-3,-1) {$\ast_3$};
        \node[anchor = west] at (3,-1) {$\ast_4$};
        \coordinate (L) at (-2,-1){};
        \coordinate (R) at (2,-1){};
        \coordinate (LL) at (-3,-1){};
        \coordinate (RR) at (3,-1){};
    \end{scope}

    \draw[yellow, line width = 0.2cm, rounded corners] (top-LL) -- (top-L) -- ($(mid-L)+(0,0.1)$) -- ($(mid-R)+(0,0.1)$) -- (bot-R) -- (bot-RR);

    \draw[cyan!50!white, line width = 0.2cm, rounded corners] (bot-LL) -- (bot-L) -- ($(mid-L)+(0,-0.1)$) -- ($(mid-R)+(0,-0.1)$) -- (top-R) -- (top-RR);
        
    \draw[->, shorten <=0.1cm, shorten >=0.1cm] (top-L) -- (mid-L);
    \draw[->, shorten <=0.1cm, shorten >=0.1cm] (top-R) -- (mid-R);
    \draw[->, shorten <=0.1cm, shorten >=0.1cm] (mid-L) -- (bot-L);
    \draw[->, shorten <=0.1cm, shorten >=0.1cm] (mid-R) -- (bot-R);

    \draw[snake it] (mid-L) -- (mid-R);
    \draw[snake it] (top-L) -- (top-LL);
    \draw[snake it] (top-R) -- (top-RR);
    \draw[snake it] (bot-L) -- (bot-LL);
    \draw[snake it] (bot-R) -- (bot-RR);
\end{tikzpicture}
            \caption{Resolution of the reaching.}
        \end{subfigure}
        \caption{Resolution of a non-intersecting auto-reaching. Some of the $\ast_i$ may be identified depending on the reaching.}
        \label{fig:resolution_auto}
    \end{figure}
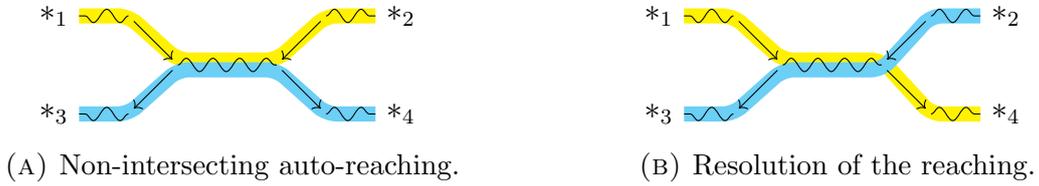
\end{definition}

\begin{example}
    We continue working with the algebra from \Cref{ex:reaching_resolution_string_band}. Let $C = \drawStringOrBand{0}{b,a,b,a,d,c,d,c}{8}$. It reaches for itself at $e_2$. The resolution is $\drawStringOrBand{0}{b,a,d,c,b,a,d,c}{8}$.
\end{example}

\begin{definition}
\label[definition]{def:intersecting_auto_reaching_resolution}
    Let $D$ be a string or band and $L = D[i,j] = D[i',j']$ an intersecting auto-reaching. We define its \df{resolution}, denoted by $\sR^{a}_L$, as the multi-set of the string or band $D_{\leq j}D_{\geq j'}$ and the band $D[j,j']$. 
\end{definition}

\begin{remark}
    The resolution of an intersecting auto-reaching takes out one of the periodicities described in \Cref{lem:intersecting_auto_reaching_periodicity}.
\end{remark}

\begin{example}
    We consider \Cref{ex:reaching_resolution_string_band}. Let $C = \drawStringOrBand{0}{b,a,b,a,d,c,b,a,d,c,d,c}{12}$. It reaches for itself at $\drawStringOrBand{0}{b,a,d,c}{4}$ which intersects in the middle of the string at $e_2$. The resolution is the string $\drawStringOrBand{0}{b,a,b,a,d,c,d,c}{8}$ and the band $\drawStringOrBand{1}{b,a,d,c}{4}$.
\end{example}

\begin{samepage}
\begin{definition}
\label[definition]{def:resolution_multi_sets}
    Let $\sD = \{D_i^{\times q_i}\}_i$ be a multi-set of bands and strings. 
    \begin{itemize}
        \item We say that $\sD$ contains an (auto-)reaching $L$ if there exist some $j$, $\ell$ and an (auto-)reaching $L$ from $D_j$ to $D_\ell$.
        \item Let $L$ be a reaching from $D_j$ to $D_\ell$ for some $(j \neq \ell)$ or $(j = \ell \text{ with } q_j>1)$. We define the resolution of $L$ in $\sD$ as $ \sR_L(\sD) := \{D_i^{\times (q_i-\delta_{ij}-\delta_{i\ell})}\}_i \bigsqcup \sR_L(D_j,D_\ell)$.
        \item Let $L$ be an auto-reaching of $D_j$ for some $j$. We define the auto-resolution of $L$ in $D$ as $ \sR^a_L(\sD) := \{D_i^{\times (q_i-\delta_{ij})}\}_i \bigsqcup \sR^a_L(D_j)$.
    \end{itemize}
\end{definition}
\end{samepage}

Resolving reachings gives rise to degenerations. The aim of this section is to prove the following.

\begin{theorem}
\label{thm::deg_and_intersection}
    Let $\sD = \{D_i^{\times q_i}\}_i$ be a multi-set of bands and strings containing an (auto-)reaching $L$. Then $ \sR^{(a)}_L(\sD) \leqdegfam \sD$. 
\end{theorem}

\begin{remark}
\label[remark]{rmk:attractions}
    In \cite{PPP_kissing}, reachings appear alongside the notion of attractions. Formally, an attraction can also be resolved. However, we do not consider them here since the dimension vector of the resulting strings and bands will differ from the original strings and bands considered. 
\end{remark}

\begin{remark}
\label[remark]{rmk:resolution_string_alg}
    In string algebras, the strings and bands obtained through a resolution might contain zero relations, i.e. are not truly strings and bands. In this case we say that the resolution does not exist. However, if the resolution exists, \Cref{thm::deg_and_intersection} holds.
\end{remark}

We will prove \Cref{thm::deg_and_intersection} in three parts: first for reachings, then for non-intersecting auto-reachings and finally for intersecting auto-reachings.

\begin{lemma}
\label[lemma]{lem:reaching_resol}
    Let $\sD = \{ D, D'\}$ be a multi-set of strings or bands with a reaching $L$ between $D$ and $D'$. Then $\sR_L \leqdegfam \sD$.
\end{lemma}

\begin{proof}
    We consider the case where $D = B$ and $D' = B'$ are two bands and $\ell(L) < \ell(B), \ell(B')$. In particular, we want to avoid $L_\text{red}$ for now. Then $\sR_L = \wB$ where $\wB$ is also a band. Due to \Cref{lem:family_closure_quasi_length_1}, it is enough to show that 
    \begin{equation*}
        \forall \lambda, \lambda' \in K^*, M(B,\lambda,1) \oplus M(B',\lambda',1) \in \overline{\sO_{\sR_L}}.
    \end{equation*}
    Similarly to the proof of \Cref{lem:deleting_arrow}, we will ultimately construct a continuous morphism between varieties $\phi: K \rightarrow \mod(A, \dd)$ where $\dd$ is the dimension vector of $\sD$ such that $\phi(0) = M(B,\lambda,1) \oplus M(B',\lambda',1)$ and $\phi(t) \in \sO_{\sR_L}$ for all $t \neq 0$. To describe $\phi$ in an efficient manner, we first build some new maps. 
    
    We call $M = M(B,\lambda,1)$ and $M' = M(B',\lambda',1)$. For $ \wlambda \in K^*$, we define the module $\wM=M(\wB, \wlambda, 1)$. Moreover, we denote by $m = \ell(B)$, $m' = \ell(B')$ and $\wm = \ell(\wB)$. Recall \Cref{constr::bands}. To build $M$, we considered a quiver morphism $f \colon \wA_{m-1} \rightarrow Q$ which induces a functor $F_f \colon \mod(\wA_{m-1}) \rightarrow \mod(A)$. For $M'$, we use $f' \colon \wA_{m'-1} \rightarrow Q$ and $F_{f'} \colon \mod(\wA_{m'-1}) \rightarrow \mod(A)$. Similarly, we denote the quiver morphism involved with $\wB$ by $\widetilde{f}$. Clearly, the functor 
    \begin{equation*}
         F_f \times F_{f'} \colon \mod(\wA_{m-1}) \times \mod(\wA_{m'-1}) \rightarrow \mod(A)
    \end{equation*}
    can be used to obtain $M \oplus M'$. 
    
    Let $\ell$, resp. $\ell'$, be the string in $\wA_{m-1}$, resp. $\wA_{m'-1}$, such that $f$, resp. $f'$, sends it isomorphically onto $L$. Using this, we identify $\ell$, $L$ and $\ell'$. We consider the coequalizer in the category of quivers, denoted $\wA_L$, of the inverse maps 
    \[\begin{tikzcd}
	L & {\wA_{m-1} \sqcup \wA_{m'-1}} & {\wA_L}.
	\arrow[shift right, from=1-1, to=1-2]
	\arrow[shift left, from=1-1, to=1-2]
	\arrow[from=1-2, to=1-3]
\end{tikzcd}\]
    This essentially glues $\ell$ and $\ell'$ together. The quiver $\wA_L$ has the form
    \[\begin{tikzcd}[sep = small]
    	\bullet &&&& \bullet \\
    	& \bullet && \bullet \\
    	\bullet &&&& \bullet
    	\arrow[curve={height=-12pt}, squiggly, no head, from=1-1, to=1-5]
    	\arrow[""{name=0, anchor=center, inner sep=0}, from=1-1, to=2-2]
    	\arrow[""{name=1, anchor=center, inner sep=0}, from=1-5, to=2-4]
    	\arrow[squiggly, no head, from=2-2, to=2-4]
    	\arrow[""{name=2, anchor=center, inner sep=0}, from=2-2, to=3-1]
    	\arrow[""{name=3, anchor=center, inner sep=0}, from=2-4, to=3-5]
    	\arrow[curve={height=12pt}, squiggly, no head, from=3-1, to=3-5]
    	\arrow[curve={height=12pt}, shorten <=5pt, shorten >=5pt, dotted, no head, from=0, to=2]
    	\arrow[curve={height=-12pt}, shorten <=5pt, shorten >=5pt, dotted, no head, from=1, to=3]
    \end{tikzcd}\]
    where the upper half stems from $\wA_{m'-1}$ and the lower part from $\wA_{m-1}$. We denote the gentle algebra with the dotted zero relations by $H$. Let $\dd_H$ be the dimension vector of $H$ with dimension $2$ on the image of $L$ and $1$ everywhere else.
    
    Similarly, we can consider the coequalizer of the two embeddings of $L$ into $\wA_{\wm-1}$. It is also $\wA_L$. Moreover, the universal property applied to $f \sqcup f'$, resp. $\widetilde{f}$ gives rise to the same quiver morphism $g \colon \wA_L \rightarrow Q$:
    \[\begin{tikzcd}
    	& {\wA_{m-1} \sqcup \wA_{m'-1}} \\
    	L & {\wA_L} & Q. \\
    	& {\wA_{\wm-1}}
    	\arrow[from=1-2, to=2-2]
    	\arrow["{f \sqcup f'}", from=1-2, to=2-3]
    	\arrow[shift left, from=2-1, to=1-2]
    	\arrow[shift right, from=2-1, to=1-2]
    	\arrow[shift right, from=2-1, to=3-2]
    	\arrow[shift left, from=2-1, to=3-2]
    	\arrow["g"', from=2-2, to=2-3]
    	\arrow["{\widetilde{f}}"', from=3-2, to=2-3]
    	\arrow[from=3-2, to=2-2]
    \end{tikzcd}\]
    In particular, the quiver morphism $g$ induces a functor $F_g$ between the module categories of $H$ and $A$ whose image contains $M \oplus M'$ and $\wM$. So $M \oplus M' = F_g (N)$ and $\wM = F_g(\wN)$ where $N,\wN \in \mod(H)$ are
    \[N = \begin{tikzcd}[sep = small]
    	K &&&& K \\
    	& K^2 && K^2 \\
    	K &&&& K,
    	\arrow["\id", curve={height=-12pt}, squiggly, no head, from=1-1, to=1-5]
    	\arrow["\pmat{\lambda' \\ 0}", ""{name=0, anchor=center, inner sep=0}, from=1-1, to=2-2]
    	\arrow["\pmat{1 \\ 0}"', ""{name=1, anchor=center, inner sep=0}, from=1-5, to=2-4]
    	\arrow["\id", squiggly, no head, from=2-2, to=2-4]
    	\arrow["\pmat{0 \amsamp \lambda}", ""{name=2, anchor=center, inner sep=0}, from=2-2, to=3-1]
    	\arrow["\pmat{0 \amsamp 1}"', ""{name=3, anchor=center, inner sep=0}, from=2-4, to=3-5]
    	\arrow["\id",curve={height=12pt}, squiggly, no head, from=3-1, to=3-5]
    	\arrow[curve={height=12pt}, shorten <=5pt, shorten >=5pt, dotted, no head, from=0, to=2]
    	\arrow[curve={height=-12pt}, shorten <=5pt, shorten >=5pt, dotted, no head, from=1, to=3]
    \end{tikzcd} \qquad \qquad \wN = \begin{tikzcd}[sep = small]
    	K &&&& K \\
    	& K^2 && K^2 \\
    	K &&&& K.
    	\arrow["\id", curve={height=-12pt}, squiggly, no head, from=1-1, to=1-5]
    	\arrow["\pmat{\wlambda \\ 0}" {yshift=-4pt}, ""{name=0, anchor=center, inner sep=0}, from=1-1, to=2-2]
    	\arrow["\pmat{0 \\ 1}"', ""{name=1, anchor=center, inner sep=0}, from=1-5, to=2-4]
    	\arrow["\id", squiggly, no head, from=2-2, to=2-4]
    	\arrow["\pmat{0 \amsamp 1}", ""{name=2, anchor=center, inner sep=0}, from=2-2, to=3-1]
    	\arrow["\pmat{1 \amsamp 0}"', ""{name=3, anchor=center, inner sep=0}, from=2-4, to=3-5]
    	\arrow["\id",curve={height=12pt}, squiggly, no head, from=3-1, to=3-5]
    	\arrow[curve={height=12pt}, shorten <=5pt, shorten >=5pt, dotted, no head, from=0, to=2]
    	\arrow[curve={height=-12pt}, shorten <=5pt, shorten >=5pt, dotted, no head, from=1, to=3]
    \end{tikzcd}\]
    We see that these modules are identical except for the four morphisms to and from $L$. So we consider the algebra $X$ given by the following quiver and relations and let $\dd_X = \pmat{1 && 1\\ &2& \\ 1&&1}$ be a dimension vector of $X$.
    \[\begin{tikzcd}
    	\bullet && \bullet \\
    	& \bullet \\
    	\bullet && \bullet
    	\arrow[""{name=0, anchor=center, inner sep=0}, from=1-1, to=2-2]
    	\arrow[""{name=1, anchor=center, inner sep=0}, from=1-3, to=2-2]
    	\arrow[""{name=2, anchor=center, inner sep=0}, from=2-2, to=3-1]
    	\arrow[""{name=3, anchor=center, inner sep=0}, from=2-2, to=3-3]
    	\arrow[curve={height=12pt}, shorten <=5pt, shorten >=5pt, dotted, no head, from=0, to=2]
    	\arrow[curve={height=-12pt}, shorten <=7pt, shorten >=7pt, dotted, no head, from=1, to=3]
    \end{tikzcd}\]
    We define a continuous morphism of varieties $\psi \colon \mod(X,\dd_X) \rightarrow \mod(H,\dd_H)$ by
    \[\begin{tikzcd}[sep = small]
    	K && K \\
    	& K^2 \\
    	K && K
    	\arrow["M_a", ""{name=0, anchor=center, inner sep=0}, from=1-1, to=2-2]
    	\arrow["M_b"', ""{name=1, anchor=center, inner sep=0}, from=1-3, to=2-2]
    	\arrow["M_c", ""{name=2, anchor=center, inner sep=0}, from=2-2, to=3-1]
    	\arrow["M_d"', ""{name=3, anchor=center, inner sep=0}, from=2-2, to=3-3]
    	\arrow[curve={height=12pt}, shorten <=5pt, shorten >=5pt, dotted, no head, from=0, to=2]
    	\arrow[curve={height=-12pt}, shorten <=7pt, shorten >=7pt, dotted, no head, from=1, to=3]
    \end{tikzcd} \longmapsto \begin{tikzcd}[sep = small]
    	K &&&& K \\
    	& K^2 && K^2 \\
    	K &&&& K.
    	\arrow["\id", curve={height=-12pt}, squiggly, no head, from=1-1, to=1-5]
    	\arrow["M_a", ""{name=0, anchor=center, inner sep=0}, from=1-1, to=2-2]
    	\arrow["M_b"', ""{name=1, anchor=center, inner sep=0}, from=1-5, to=2-4]
    	\arrow["\id", squiggly, no head, from=2-2, to=2-4]
    	\arrow["M_c", ""{name=2, anchor=center, inner sep=0}, from=2-2, to=3-1]
    	\arrow["M_d"', ""{name=3, anchor=center, inner sep=0}, from=2-4, to=3-5]
    	\arrow["\id",curve={height=12pt}, squiggly, no head, from=3-1, to=3-5]
    	\arrow[curve={height=12pt}, shorten <=5pt, shorten >=5pt, dotted, no head, from=0, to=2]
    	\arrow[curve={height=-12pt}, shorten <=5pt, shorten >=5pt, dotted, no head, from=1, to=3]
    \end{tikzcd}\]
    We define $M_{\lambda,\lambda'}\in \mod(X,\dd_X)$ by $M_a = \pmat{\lambda' \\ 0}$, $M_b = \pmat{1 \\ 0}$, $M_c = \pmat{0 & \lambda}$ and $M_d = \pmat{0 & 1}$. By a slight abuse of notation, we denote by $F_g$ the morphism between module varieties induced by $F_g$ as well. Then $F_g(\psi(M_{\lambda, \lambda'})) = M \oplus M'$. We will build a morphism of varieties $\varphi \colon K \rightarrow \mod(X, \dd_X)$ such that $\varphi(0) = M_{\lambda,\lambda'}$ and $F_g(\psi(\varphi(t))) \in \sO_{\sR_L}$. This will imply that $M\oplus M' \in \overline{\sO_{\sR_L}}$. For all $t\in K$, we define $\varphi(t)$ as
    \begin{equation}
    \label{quiver:morphism_variety}
        \begin{tikzcd}
            K && K \\
            & K^2 \\
            K && K.
            \arrow["\pmat{\lambda' \\ 0}"', from=1-1, to=2-2]
            \arrow["\pmat{1-t \\ -t}"', from=1-3, to=2-2]
            \arrow["\pmat{0 \amsamp \lambda}", from=2-2, to=3-1]
            \arrow["\pmat{t \amsamp 1-t}", from=2-2, to=3-3]
        \end{tikzcd}
    \end{equation}
    For $t=0$, this is indeed $M_{\lambda, \lambda'}$. If $t \neq 0$, we can apply a base change at $K^2$ to $\left\{ \pmat{1 \\ 0}, \pmat{1-t \\ -t} \right\}$ to show 
    \[\varphi(t) \cong \begin{tikzcd}
        K && K \\
        & K^2 \\
        K && K.
        \arrow["\pmat{\lambda' \\ 0}"', from=1-1, to=2-2]
        \arrow["\pmat{0 \\ 1}"', from=1-3, to=2-2]
        \arrow["\pmat{0 \amsamp -\lambda t}", from=2-2, to=3-1]
        \arrow["\pmat{t \amsamp 0}", from=2-2, to=3-3]
    \end{tikzcd}\]
    This isomorphism is preserved by $F_g \circ \psi$. Moreover, the above module gets sent to the family of $\sR_L$. In particular, so does $\varphi(t)$.
    
    If $\ell(L) \geq m$, or $m'$, we work with up to two $L_\text{red}$. So we cannot take the coequalizer for the embedding of $L$, but the one for the embedding of the endpoints of $L$, generating the following quiver.
    \[\begin{tikzcd}[sep = small]
    	\bullet &&&& \bullet \\
    	& \bullet && \bullet \\
    	\bullet &&&& \bullet
    	\arrow[curve={height=-12pt}, squiggly, no head, from=1-1, to=1-5]
    	\arrow[""{name=0, anchor=center, inner sep=0}, from=1-1, to=2-2]
    	\arrow[""{name=1, anchor=center, inner sep=0}, from=1-5, to=2-4]
    	\arrow[curve={height=-12pt}, shift right, squiggly, no head, from=2-2, to=2-4]
    	\arrow[curve={height=12pt}, shift left, squiggly, no head, from=2-2, to=2-4]
    	\arrow[""{name=2, anchor=center, inner sep=0}, from=2-2, to=3-1]
    	\arrow[""{name=3, anchor=center, inner sep=0}, from=2-4, to=3-5]
    	\arrow[curve={height=12pt}, squiggly, no head, from=3-1, to=3-5]
    	\arrow[curve={height=12pt}, shorten <=5pt, shorten >=5pt, dotted, no head, from=0, to=2]
    	\arrow[curve={height=-12pt}, shorten <=5pt, shorten >=5pt, dotted, no head, from=1, to=3]
    \end{tikzcd}\]
    The reduction to $X$ remains possible and thus the proof can be completed. For the case where one or two strings are involved in the reaching, the above proof can be adapted by replacing the appropriate $\wA_{m-1}$ with $A_{m+1}$. 
\end{proof}

\begin{lemma}
\label[lemma]{lem:non_intersecting_auto_reaching_resol}
    Let $D$ be a  string or band with a non-intersecting auto-reaching $L$ on $D$. Then $\sR^{a}_L \leqdegfam D$.
\end{lemma}

\begin{proof}
    The proof of \Cref{lem:reaching_resol} can be adapted to this case. Note that it is possible to have reachings where the two substrings are only one arrow apart. In this case, we compose the two morphisms we would usually assign to the two arrows in (\ref{quiver:morphism_variety}), for example
    \[\begin{tikzcd}
    	K \\
    	& {K^2} \\
    	K
    	\arrow["\pmat{\lambda' \\ 0}", from=1-1, to=2-2]
    	\arrow["\pmat{1-t \\ -t} \circ \pmat{t \amsamp 1-t}", from=2-2, to=2-2, loop, in=325, out=35, distance=10mm]
    	\arrow["\pmat{0 \amsamp \lambda}", from=2-2, to=3-1]
    \end{tikzcd}\]
    where the two right arrows got replaced by a loop.
\end{proof}

\begin{lemma}
\label[lemma]{lem:intersecting_auto_reaching_resol}
    Let $D$ be a  string or band with an intersecting auto-reaching $L$ on $D$. Then $\sR^{a}_L \leqdegfam D$.
\end{lemma}

\begin{proof}
    Again, we may follow the proof of \Cref{lem:reaching_resol}. However, some details change. Let us consider the case of a band $B$ with an intersecting auto-reaching $L$. By \Cref{lem:intersecting_auto_reaching_periodicity}, we know that $L = \alpha^s$ for some string $\alpha$ and periodicity $s$ such that $L$ overlaps on $\alpha^{s-1}$. Due to this, forming the coequalizer of the two copies of $L$ on $\wA_{m-1}$ where $m = \ell(B)$ is equivalent with identifying all these $s+1$ copies of $\alpha$. In particular, $\wA_L$ is a copy of $\wA_{m - (s+1) k}$ and $\wA_k$ glued together at one vertex, see \Cref{fig:gluing_intersecting_string}. Here, $k = \ell(\alpha)$.
    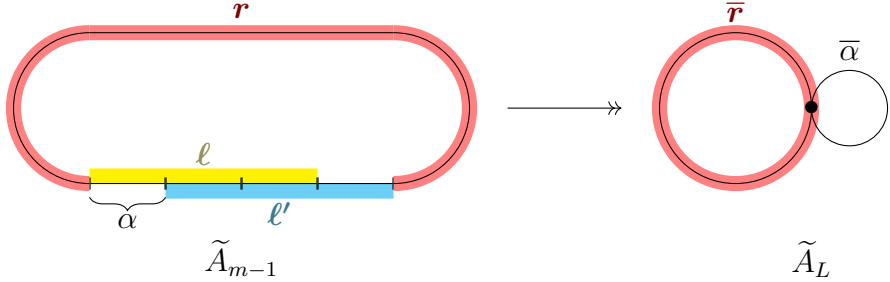
\begin{figure}[H]
        \centering
        \begin{tikzpicture}[scale = 0.5]
    
    \draw[red!50!white, line width = 0.2cm] (0,4) -- (8,4);
    \centerarc[red!50!white, line width = 0.2cm](0,2)(90:270:2);
    \centerarc[red!50!white, line width = 0.2cm](8,2)(-90:90:2);
    \draw[yellow, line width = 0.2cm] ($(0,0)+(0,0.2cm)$) -- ($(6,0)+(0,0.2cm)$);
    \draw[cyan!50!white, line width = 0.2cm] ($(2,0)+(0,-0.2cm)$) -- ($(8,0)+(0,-0.2cm)$); 
    
    \node[] at (4,-2) {$\wA_{m-1}$};
    \node[anchor = south, color = red!50!black] at ($(4,4)+(0,0.1cm)$) {$\boldsymbol{r}$};
    \node[anchor = south, color = yellow!50!black] at ($(3,0)+(0,0.2cm)$) {$\boldsymbol{\ell}$};
    \node[anchor = north, color = cyan!50!black] at ($(5,0)+(0,-0.2cm)$) {$\boldsymbol{\ell'}$};
    \draw [decorate,decoration={brace,amplitude=5pt,mirror,raise=1ex}]
  (0,0) -- (2,0) node[midway,yshift=-1.2em]{$\alpha$};
    
    \draw[|-|] (0,0) -- (2,0);  
    \draw[|-|] (2,0) -- (4,0);  
    \draw[|-|] (4,0) -- (6,0);  
    \draw[|-|] (6,0) -- (8,0); 
    \draw[] (0,4) -- (8,4);
    
    \centerarc[](0,2)(90:270:2);
    \centerarc[](8,2)(-90:90:2);
    
    \draw[->>] (11,2) -- (14,2);
    
    \centerarc[red!50!white, line width = 0.2cm](17,2)(0:360:2);
    
    \node[] at (19,2) {$\bullet$};
    \centerarc[](17,2)(0:360:2);
    \centerarc[](20,2)(0:360:1);
    
    \node[anchor = south, color = red!50!black] at (17,4) {$\boldsymbol{\overline{r}}$};
    \node[anchor = south] at (20,3) {$\overline{\alpha}$};
    \node[] at (19,-2) {$\wA_L$};
    
\end{tikzpicture}
        \caption{Gluing two intersecting strings in $\wA_{m-1}$.}
        \label{fig:gluing_intersecting_string}
    \end{figure}
    Again, we can construct an algebra $H$ with two additional relations. The module $M(B, \lambda, 1)$ will then be the image of an $\wA_L$ module $M$ of the following form
    \begin{equation*}
        \begin{tikzcd}[row sep=large, column sep=large,
            execute at end picture={
              \begin{pgfonlayer}{background}
                    \coordinate (a) at (A);
                    \coordinate (b) at (B);
                    \coordinate (c) at (C);
                    \draw[red!50!white, line width = 0.2cm, rounded corners] (a) ..controls +(-2,-5).. (c) -- (b) -- cycle;
                \end{pgfonlayer}
            }%
            ]%
        	|[alias=A]|{K} && K^{s+1} \\
        	& |[alias=B]|{K^{s+2}} \\
        	K^{s+1} && |[alias=C]|{K}
        	\arrow[""{name=0, anchor=center, inner sep=0}, "\pmat{\lambda \\ 0 \\ \Vdots \\ \Block[draw=blue]{2-1}{} 0 \\ 0}", from=1-1, to=2-2]
        	\arrow[""{name=1, anchor=center, inner sep=0}, "\pmat{1 \amsamp \amsamp \\ \amsamp \Ddots \amsamp \\ \amsamp \amsamp \Block[draw=blue]{2-1}{} 1 \\ 0 \amsamp \cdots \amsamp 0}"', from=1-3, to=2-2]
        	\arrow[""{name=2, anchor=center, inner sep=0}, "\pmat{0 \amsamp 1 \amsamp \amsamp \amsamp \\ \Vdots \amsamp \amsamp \Ddots \amsamp \amsamp \\  \amsamp \amsamp \amsamp 1 \amsamp \\ 0 \amsamp \amsamp \amsamp \Block[draw=blue]{1-2}{} 0 \amsamp 1}"', from=2-2, to=3-1]
        	\arrow[""{name=3, anchor=center, inner sep=0}, "\pmat{0 \amsamp \Cdots \amsamp \Block[draw=blue]{1-2}{} 0 \amsamp 1}", from=2-2, to=3-3]
        	\arrow[no head, from=1-3, to=3-1, to path={..controls +(2,-5).. (\tikztotarget)}]
	        \arrow[no head, from=1-1, to=3-3, to path={..controls +(-2,-5).. (\tikztotarget)}]
        \end{tikzcd}
    \end{equation*}
    with identities on the unlabeled parts. The coloured circle is the non glued part of $\wA_L$. The resolution of the reaching is defined by taking out one of the copies of $\alpha$ in $B$ and considering it as a second band. We will morally remove the last copy of $\alpha$. The module associated to the resolution is 
    \begin{equation*}
        \begin{tikzcd}[row sep=large, column sep=large,
            execute at end picture={
              \begin{pgfonlayer}{background}
                    \coordinate (a) at (A);
                    \coordinate (b) at (B);
                    \coordinate (c) at (C);
                    \draw[red!50!white, line width = 0.2cm, rounded corners] (a) ..controls +(-2,-5).. (c) -- (b) -- cycle;
                \end{pgfonlayer}
            }%
            ]%
        	|[alias=A]|{K} && K^{s+1} \\
        	& |[alias=B]|{K^{s+2}} \\
        	K^{s+1} && |[alias=C]|{K}.
        	\arrow[""{name=0, anchor=center, inner sep=0}, "\pmat{\lambda \\ 0 \\ \Vdots \\ \Block[draw=blue]{2-1}{} 0 \\ 0}", from=1-1, to=2-2]
        	\arrow[""{name=1, anchor=center, inner sep=0}, "\pmat{1 \amsamp \amsamp \\ \amsamp \Ddots \amsamp \\ \amsamp \amsamp \Block[draw=blue]{2-1}{} 0 \\ 0 \amsamp \cdots \amsamp 1}"', from=1-3, to=2-2]
        	\arrow[""{name=2, anchor=center, inner sep=0}, "\pmat{0 \amsamp 1 \amsamp \amsamp \amsamp \\ \Vdots \amsamp \amsamp \Ddots \amsamp \amsamp \\  \amsamp \amsamp \amsamp 1 \amsamp \\ 0 \amsamp \amsamp \amsamp \Block[draw=blue]{1-2}{} 0 \amsamp 1}"', from=2-2, to=3-1]
        	\arrow[""{name=3, anchor=center, inner sep=0}, "\pmat{0 \amsamp \Cdots \amsamp \Block[draw=blue]{1-2}{} 1 \amsamp 0}", from=2-2, to=3-3]
        	\arrow[no head, from=1-3, to=3-1, to path={..controls +(2,-5).. (\tikztotarget)}]
	        \arrow[no head, from=1-1, to=3-3, to path={..controls +(-2,-5).. (\tikztotarget)}]
        \end{tikzcd}
    \end{equation*}
    Note that only the action on the last two basis vectors of $K^{s+2}$ were affected, see the marked coloured boxes in both representations. So we can reduce to the algebra $X$ by replacing the morphism $\psi$ by one defined via 
    \[\begin{tikzcd}[sep = small]
    	K && K \\
    	& K^2 \\
    	K && K
    	\arrow["\pmat{{a_1} \\ {a_2} }", ""{name=0, anchor=center, inner sep=0}, from=1-1, to=2-2]
    	\arrow["\pmat{{b_1} \\ {b_2} }"', ""{name=1, anchor=center, inner sep=0}, from=1-3, to=2-2]
    	\arrow["\pmat{{c_1} \amsamp {c_2} }"', ""{name=2, anchor=center, inner sep=0}, from=2-2, to=3-1]
    	\arrow["\pmat{{d_1} \amsamp {d_2} }", ""{name=3, anchor=center, inner sep=0}, from=2-2, to=3-3]
    \end{tikzcd} \longmapsto \begin{tikzcd}[row sep=large, column sep=large,]
    	{K} && K^{s+1} \\
        	& {K^{s+2}} \\
        	K^{s+1} && {K}.
        	\arrow[""{name=0, anchor=center, inner sep=0}, "\pmat{\lambda \\ 0 \\ \Vdots \\ \Block[draw=blue]{2-1}{} {a_1} \\ {a_2} }", from=1-1, to=2-2]
        	\arrow[""{name=1, anchor=center, inner sep=0}, "\pmat{1 \amsamp \amsamp \\ \amsamp \Ddots \amsamp \\ \amsamp \amsamp \Block[draw=blue]{2-1}{} {b_1} \\ \amsamp \cdots \amsamp b_2}"', from=1-3, to=2-2]
        	\arrow[""{name=2, anchor=center, inner sep=0}, "\pmat{0 \amsamp 1 \amsamp \amsamp \amsamp \\ \Vdots \amsamp \amsamp \Ddots \amsamp \amsamp \\  \amsamp \amsamp \amsamp 1 \amsamp \\ 0 \amsamp \amsamp \amsamp \Block[draw=blue]{1-2}{} {c_1} \amsamp {c_2}}"', from=2-2, to=3-1]
        	\arrow[""{name=3, anchor=center, inner sep=0}, "\pmat{0 \amsamp \Cdots \amsamp \Block[draw=blue]{1-2}{} {d_1} \amsamp {d_2}}", from=2-2, to=3-3]
        	\arrow[no head, from=1-3, to=3-1, to path={..controls +(3,-5).. (\tikztotarget)}]
	        \arrow[no head, from=1-1, to=3-3, to path={..controls +(-3,-5).. (\tikztotarget)}]
    \end{tikzcd}\]
    and complete the proof as before. 
    The same can be done in the case of a string.
\end{proof}
    
Applying \Cref{lem:disjoint_union} to \Cref{lem:reaching_resol,,lem:non_intersecting_auto_reaching_resol,,lem:intersecting_auto_reaching_resol} completes the proof of \Cref{thm::deg_and_intersection}.

\biblio{}

\section{Examples}
\label{sect:examples}
We study some examples of degenerations for the algebras in \Cref{ex:band,,ex:reaching_resolution_string_band}. In particular, we give the complete partial order of degenerations for some dimension vectors. 

We start by a short example on degenerations arising from reachings. Similarly to what we observed for deletions of arrows, degenerations of strings and bands arising from reachings do not usually implicate degeneration between individual modules.
\begin{example}
\label[example]{ex:band_deg_vs_mod_deg}
    Let $A$ be the algebra given by quiver and relations
    \[\begin{tikzcd}[sep = small]
    	1 && 2 \\
    	& 3 \\
    	4 && 5.
    	\arrow[""{name=0, anchor=center, inner sep=0}, "a", from=1-1, to=2-2]
    	\arrow[""{name=1, anchor=center, inner sep=0}, "b"', from=1-3, to=2-2]
    	\arrow[""{name=2, anchor=center, inner sep=0}, "c", from=2-2, to=3-1]
    	\arrow[""{name=3, anchor=center, inner sep=0}, "d"', from=2-2, to=3-3]
    	\arrow["{e}"', curve={height=30pt}, from=1-1, to=3-1]
    	\arrow["{f}", curve={height=-30pt}, from=1-3, to=3-3]
    	\arrow[curve={height=-12pt}, shorten <=5pt, shorten >=5pt, dashed, no head, from=1, to=3]
    	\arrow[curve={height=12pt}, shorten <=5pt, shorten >=5pt, dashed, no head, from=0, to=2]
    \end{tikzcd}\]
    We consider the band
    \begin{equation*}
        B' = 
        \begin{tikzcd}[sep = small]
        	& 2 && 3 && 1 \\
        	\circled{3} && 5 && 4 && \circled{3}.
        	\arrow["a", from=1-6, to=2-7]
        	\arrow["e"', from=1-6, to=2-5]
        	\arrow["c", from=1-4, to=2-5]
        	\arrow["b"', from=1-2, to=2-1]
        	\arrow["f", from=1-2, to=2-3]
        	\arrow["d"', from=1-4, to=2-3]
        \end{tikzcd}
    \end{equation*}
    It reaches for itself at vertex $3$. The resolution of this reaching is the following band
    \begin{equation*}
        B = 
        \begin{tikzcd}[sep = small]
        	\circled{3} && 2 \\
        	& 5 && 3 && 1 \\
        	&&&& 4 && \circled{3}.
        	\arrow["a", from=2-6, to=3-7]
        	\arrow["e"', from=2-6, to=3-5]
        	\arrow["c", from=2-4, to=3-5]
        	\arrow["b", from=1-3, to=2-4]
        	\arrow["f"', from=1-3, to=2-2]
        	\arrow["d", from=1-1, to=2-2]
        \end{tikzcd}
    \end{equation*}
    As a consequence of \Cref{thm::deg_and_intersection}, $B \leqdegfam B'$. However, for any $\lambda, \lambda' \in K^*$, there are no homomorphisms in either direction between the band modules $M(B, \lambda, 1)$ and $M(B', \lambda', 1)$. By upper semicontinuity of $\hom$, it follows that there is no degeneration between these modules. This is an example of a degeneration of a family of band modules which is not implied by a classical degeneration of modules.
\end{example}

We saw in \Cref{ex:resolution_produces_non_minimal_bands} that resolutions of minimal bands do not necessarily produce minimal bands again. However, this theory can be reduced to minimal bands. By \Cref{lem:non_minimal_band_deg}, for any minimal band $B$ and any $q>1$, $B^{\times q} \leqdegfam B^q$. So any resolution of a reaching can further be reduced to minimal bands. Let $\sD = \{D_i^{\times q_i}\}_i$ be a multi-set of minimal bands and strings containing an (auto-)reaching $L$. By \Cref{thm::deg_and_intersection}, $ \sR^{(a)}_L(\sD) \leqdegfam \sD$. If $\sR^{(a)}_L(\sD)$ contains non-minimal bands, we can replace them by the appropriate amount of copies of minimal bands to obtain a multi-set of minimal-bands and strings $(\sR^{(a)}_L(\sD))_{\text{min}}$. Then $(\sR^{(a)}_L(\sD))_{\text{min}} \leqdeg \sR^{(a)}_L(\sD)$ by \Cref{lem:non_minimal_band_deg}. By restricting to minimal bands, we will say that the reaching $L$ induces the resolution $(\sR^{(a)}_L(\sD))_{\text{min}} \leqdegfam \sD$. We do this to restrict the size of all following examples.

\begin{example}
\label[example]{ex:Kronecker_dim22}
    Let us return to the Kronecker quiver, see \Cref{ex:band}. The degenerations for dimension vector $(2,2)$ form the poset in \Cref{fig:Kronecker_dim22}. All covering relations stem from either forgetting a single arrow or resolving of a reaching. The latter is marked by squiggly lines. Note that here, the degeneration order is the same as the $h$-order. We give the values of the $h$-vectors for strings $(e_1,a,b,a^-ba^-,ba^-b,ba^-ba^-)$ underneath the multi-sets.
\begin{figure}[H]
    \centering
\[\begin{tikzcd}[sep = {20mm,between origins}]
	&&& 
	{\stackquiver{
	\left\{e_1^{\times 2}, e_2^{\times 2}\right\} 
	\\ {\scriptstyle\textcolor{teal!50!blue}{(2, 2, 2, 4, 4, 4)}}}}
	&&& \\
	&&
	{\stackquiver{
	\left\{\drawStringOrBandKronecker{0}{a}{1}, e_1, e_2\right\}
	\\ {\scriptstyle\textcolor{teal!50!blue}{(1, 2, 1, 3, 2, 3)}}}}
	&&
	{\stackquiver{
	\left\{ \drawStringOrBandKronecker{0}{b}{1}, e_1, e_2\right\}
	\\ {\scriptstyle\textcolor{teal!50!blue}{(1, 1, 2, 2, 3, 3)}}}}
	\\
	&&& 
	{\stackquiver{
	\left\{\drawStringOrBandKronecker{1}{a,b}{2}, e_1, e_2\right\}
	\\ {\scriptstyle\textcolor{teal!50!blue}{(1, 1, 1, 2, 2, 3)}}}}
	\\
	{\stackquiver{
	\left\{\drawStringOrBandKronecker{0}{a}{1}, \drawStringOrBandKronecker{0}{a}{1}\right\}
	\\ {\scriptstyle\textcolor{teal!50!blue}{(0, 2, 0, 2, 0, 2)}}}}
	&& 
	{\stackquiver{
	\left\{\drawStringOrBandKronecker{0}{a,b}{2}, e_2\right\}
	\\ {\scriptstyle\textcolor{teal!50!blue}{(0, 1, 1, 2, 2, 3)}}}}
	&& 
	{\stackquiver{
	\left\{\drawStringOrBandKronecker{0}{b,a}{2}, e_1\right\} 
	\\ {\scriptstyle\textcolor{teal!50!blue}{(1, 1, 1, 2, 2, 2)}}}}
	&& 
	{\stackquiver{
	\left\{\drawStringOrBandKronecker{0}{b}{1}, \drawStringOrBandKronecker{0}{b}{1}\right\}
	\\ {\scriptstyle\textcolor{teal!50!blue}{(0, 0, 2, 0, 2, 2)}}}}
	\\
	& 
	{\stackquiver{
	\left\{\drawStringOrBandKronecker{0}{a,b,a}{3}\right\} 
	\\ {\scriptstyle\textcolor{teal!50!blue}{(0, 1, 0, 2, 0, 2)}}}}
	&& 
	{\stackquiver{
	\left\{\drawStringOrBandKronecker{0}{a}{1}, \drawStringOrBandKronecker{0}{b}{1}\right\}
	\\ {\scriptstyle\textcolor{teal!50!blue}{(0, 1, 1, 1, 1, 2)}}}}
	&& 
	{\stackquiver{
	\left\{\drawStringOrBandKronecker{0}{b,a,b}{3}\right\} 
	\\ {\scriptstyle\textcolor{teal!50!blue}{(0, 0, 1, 0, 2, 2)}}}}
	\\
	&& 
	{\stackquiver{
	\left\{\drawStringOrBandKronecker{1}{a,b}{2}, \drawStringOrBandKronecker{0}{a}{1}\right\}
	\\ {\scriptstyle\textcolor{teal!50!blue}{(0, 1, 0, 1, 0, 2)}}}}
	&& 
	{\stackquiver{
	\left\{\drawStringOrBandKronecker{1}{a,b}{2}, \drawStringOrBandKronecker{0}{b}{1}\right\}
	\\ {\scriptstyle\textcolor{teal!50!blue}{(0, 0, 1, 0, 1, 2)}}}}
	\\
	&&&
	{\stackquiver{
	\left\{\drawStringOrBandKronecker{1}{a,b}{2},\drawStringOrBandKronecker{1}{a,b}{2} \right\}
	\\ {\scriptstyle\textcolor{teal!50!blue}{(0, 0, 0, 0, 0, 2)}}}}
	\arrow[no head, from=1-4, to=2-3]
	\arrow[no head, from=1-4, to=2-5]
	\arrow[no head, from=2-3, to=3-4]
	\arrow[no head, from=2-3, to=4-1]
	\arrow[no head, from=2-5, to=3-4]
	\arrow[no head, from=2-5, to=4-7]
	\arrow[squiggly, no head, from=3-4, to=4-3]
	\arrow[squiggly, no head, from=3-4, to=4-5]
	\arrow[no head, from=4-1, to=5-2]
	\arrow[no head, from=4-3, to=5-2]
	\arrow[squiggly, no head, from=4-3, to=5-4]
	\arrow[squiggly, no head, from=4-5, to=5-4]
	\arrow[no head, from=4-5, to=5-6]
	\arrow[no head, from=4-7, to=5-6]
	\arrow[squiggly, no head, from=5-2, to=6-3]
	\arrow[no head, from=5-4, to=6-3]
	\arrow[no head, from=5-4, to=6-5]
	\arrow[squiggly, no head, from=5-6, to=6-5]
	\arrow[no head, from=6-3, to=7-4]
	\arrow[no head, from=6-5, to=7-4]
\end{tikzcd}\]
    \caption{Degeneration order of the Kronecker quiver for $\dd = (2,2)$. The squiggly lines represent the resolutions of reachings. The coloured vectors show certain entries of the $h$-vector.}
    \label{fig:Kronecker_dim22}
\end{figure}
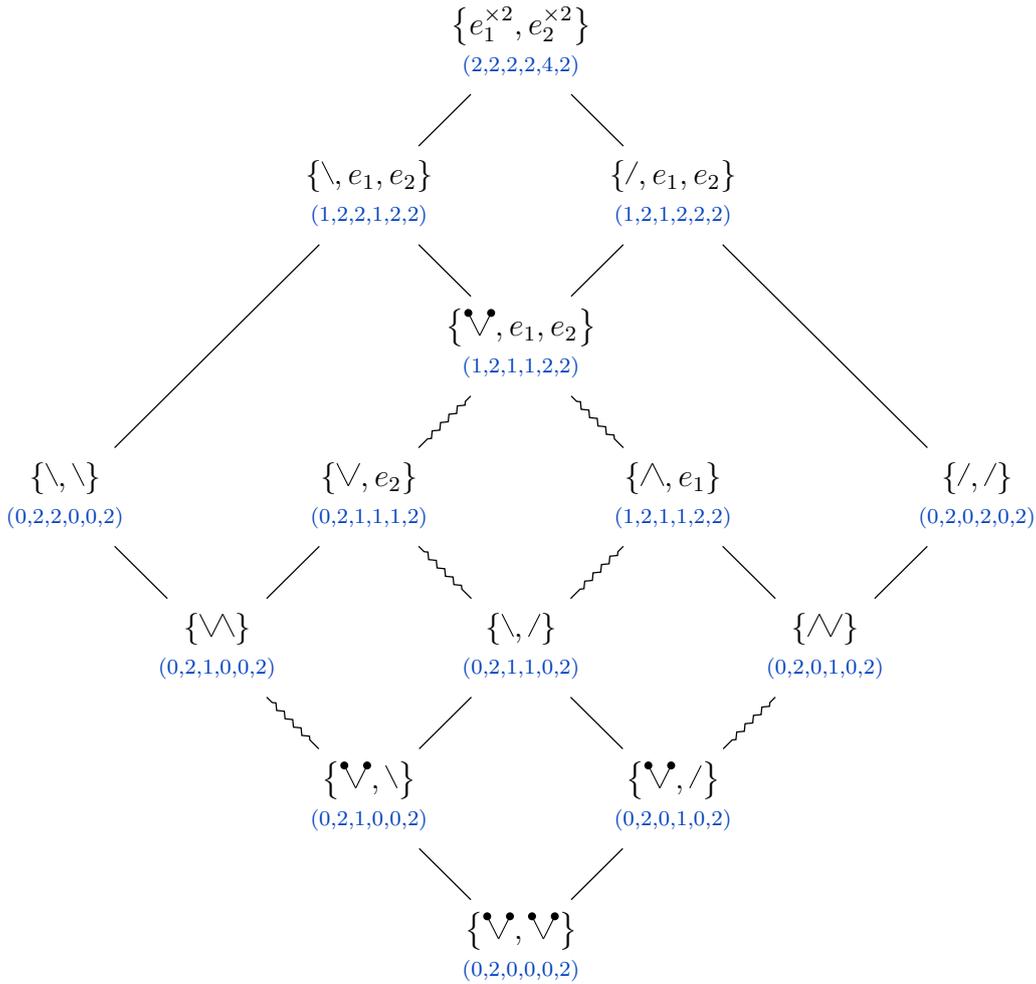
\end{example}

\begin{figure}[H]
    \centering
    \includegraphics[height = 7 cm]{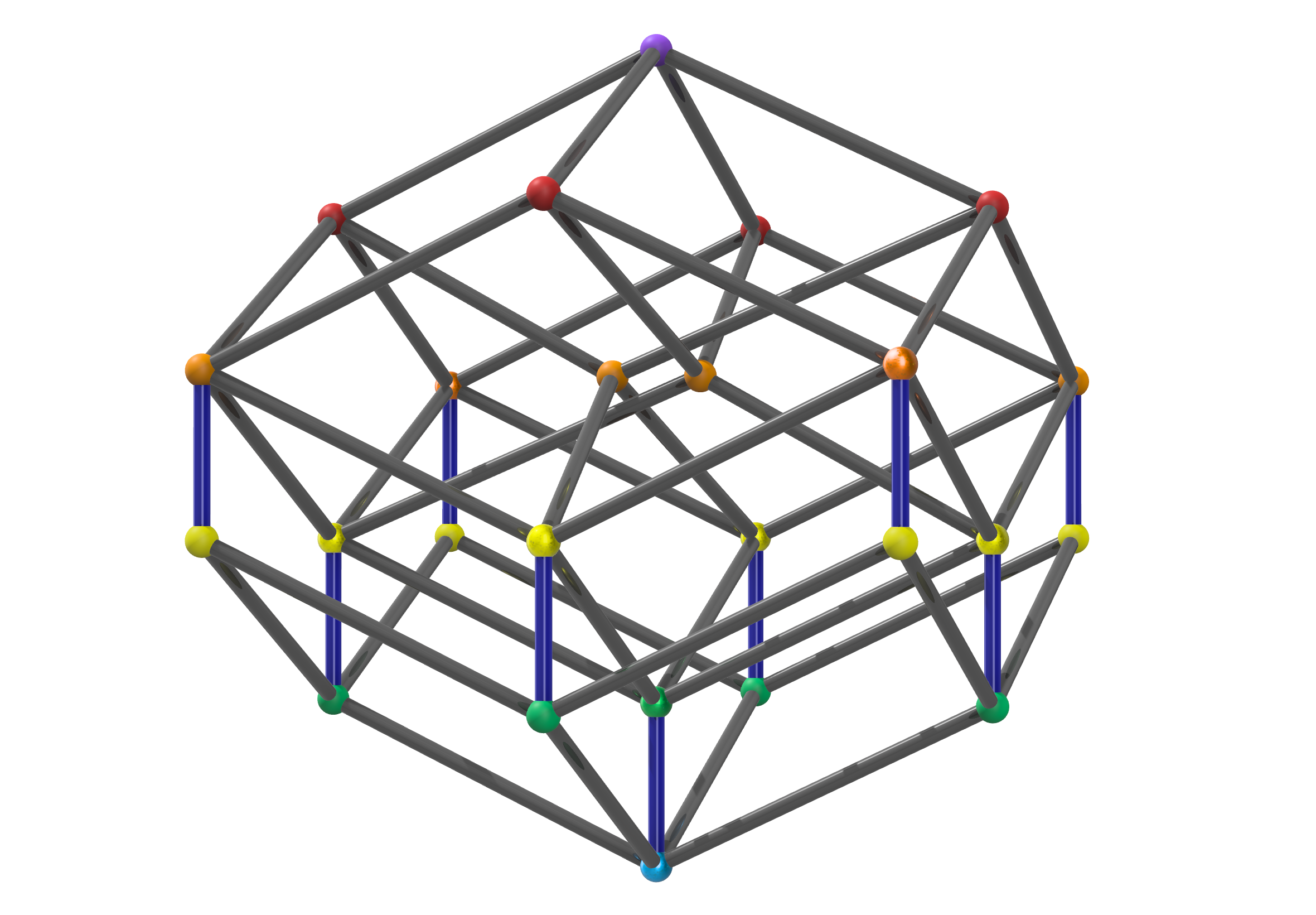}
    \caption{Poset spanned by degenerations of strings and bands for dimension vector $(1,2,1)$ in 3D. The coloured edges represent the resolutions of reachings.}
    \label{fig:deg_dim_121_3D}
\end{figure}

\begin{example}
\label[example]{ex:glued_Kronecker_dim121}
    We reconsider the algebra in \Cref{ex:reaching_resolution_string_band}. For the dimension vector $(1,2,1)$, the poset of degenerations is shown in \Cref{fig:deg_dim_121}. Again, covering relations stemming from resolutions are marked by squiggly lines, all other covering relations are obtained by deleting an arrow. Note that there are instances where deleting arrows does not produce cover relations, see the rightmost mutli-sets in the second and fourth row from the top. 
    We also note that for this example, the $h$-order and degeneration order coincide again. We give the $h$-vectors for entries $(e_1,e_2,c,d,da,cb,a^-ba^-b)$.
    This poset is not a planar graph, but it can be modelled without intersections in three dimensions, see \Cref{fig:deg_dim_121_3D}. The different ranks of the poset are differentiated by colour. Moreover, the depiction aligns with the one before; for every rank, the vertices from left to right are the multi-sets in \Cref{fig:deg_dim_121} from left to right. The resolutions are marked by coloured edges.
\end{example}

\begin{landscape}%
    \begin{figure}[p]%
        \centering%
        \[\begin{tikzcd}[column sep=tiny, ampersand replacement=\&]
	\&\&\&\&
	{\stackquiver{
	{\left\{e_1,e_2^{\times 2},e_3\right\}}
	\\ {\scriptstyle\textcolor{teal!50!blue}{(1, 2, 2, 2, 1, 1, 2)}}}}
	\\
	\& 
	{\stackquiver{
	{\left\{\drawStringOrBand{0}{d}{1}, e_1, e_2\right\}}
	\\ {\scriptstyle\textcolor{teal!50!blue}{(1, 1, 1, 2, 1, 1, 2)}}}}
	\&\& 
	{\stackquiver{
	{\left\{\drawStringOrBand{0}{a}{1}, e_2, e_3\right\}}
	\\ {\scriptstyle\textcolor{teal!50!blue}{(0, 2, 2, 2, 1, 0, 1)}}}}
	\&\& 
	{\stackquiver{
	{\left\{\drawStringOrBand{0}{c}{1}, e_1, e_2\right\}}
	\\ {\scriptstyle\textcolor{teal!50!blue}{(1, 1, 2, 1, 1, 1, 2)}}}}
	\&\& 
	{\stackquiver{
	{\left\{\drawStringOrBand{0}{b}{1}, e_2, e_3\right\}}
	\\ {\scriptstyle\textcolor{teal!50!blue}{(0, 2, 2, 2, 0, 1, 1)}}}}
	\\
	{\stackquiver{
	{\left\{\drawStringOrBand{0}{a}{1}, \drawStringOrBand{0}{d}{1}\right\}}
	\\ {\scriptstyle\textcolor{teal!50!blue}{(0, 1, 1, 2, 1, 0, 1)}}}}
	\&\& 
	{\stackquiver{
	{\left\{\drawStringOrBand{1}{d,c}{2}, e_1, e_2\right\}}
	\\ {\scriptstyle\textcolor{teal!50!blue}{(1, 1, 1, 1, 1, 1, 2)}}}}
	\& 
	{\stackquiver{
	{\left\{\drawStringOrBand{0}{b}{1}, \drawStringOrBand{0}{d}{1}\right\}}
	\\ {\scriptstyle\textcolor{teal!50!blue}{(0, 1, 1, 2, 0, 1, 1)}}}}
	\&\& 
	{\stackquiver{
	{\left\{\drawStringOrBand{0}{a}{1}, \drawStringOrBand{0}{c}{1}\right\}}
	\\ {\scriptstyle\textcolor{teal!50!blue}{(0, 1, 2, 1, 1, 0, 1)}}}}
	\& 
	{\stackquiver{
	{\left\{\drawStringOrBand{1}{b,a}{2}, e_2, e_3\right\}}
	\\ {\scriptstyle\textcolor{teal!50!blue}{(0, 2, 2, 2, 0, 0, 1)}}}}
	\&\& 
	{\stackquiver{
	{\left\{\drawStringOrBand{0}{b}{1}, \drawStringOrBand{0}{c}{1}\right\}}
	\\ {\scriptstyle\textcolor{teal!50!blue}{(0, 1, 2, 1, 0, 1, 1)}}}}
	\\
	{\stackquiver{
	{\left\{\drawStringOrBand{0}{a,d}{2}, e_2\right\}}
	\\ {\scriptstyle\textcolor{teal!50!blue}{(0, 1, 1, 2, 1, 0, 0)}}}}
	\& 
	{\stackquiver{
	{\left\{\drawStringOrBand{0}{a}{1}, \drawStringOrBand{1}{d,c}{2}\right\}}
	\\ {\scriptstyle\textcolor{teal!50!blue}{(0, 1, 1, 1, 1, 0, 1)}}}}
	\& 
	{\stackquiver{
	{\left\{\drawStringOrBand{0}{d,c}{2}, e_1\right\}}
	\\ {\scriptstyle\textcolor{teal!50!blue}{(1, 0, 1, 1, 1, 1, 2)}}}}
	\& 
	{\stackquiver{
	{\left\{\drawStringOrBand{1}{b,a}{2}, \drawStringOrBand{0}{d}{1}\right\}}
	\\ {\scriptstyle\textcolor{teal!50!blue}{(0, 1, 1, 2, 0, 0, 1)}}}}
	\&\& 
	{\stackquiver{
	{\left\{\drawStringOrBand{0}{b}{1}, \drawStringOrBand{1}{d,c}{2}\right\}}
	\\ {\scriptstyle\textcolor{teal!50!blue}{(0, 1, 1, 1, 0, 1, 1)}}}}
	\& 
	{\stackquiver{
	{\left\{\drawStringOrBand{0}{b,a}{2}, e_3\right\}}
	\\ {\scriptstyle\textcolor{teal!50!blue}{(0, 2, 2, 2, 0, 0, 0)}}}}
	\& 
	{\stackquiver{
	{\left\{\drawStringOrBand{1}{b,a}{2}, \drawStringOrBand{0}{c}{1}\right\}}
	\\ {\scriptstyle\textcolor{teal!50!blue}{(0, 1, 2, 1, 0, 0, 1)}}}}
	\& 
	{\stackquiver{
	{\left\{\drawStringOrBand{0}{c,b}{2}, e_2\right\}}
	\\ {\scriptstyle\textcolor{teal!50!blue}{(0, 1, 2, 1, 0, 1, 0)}}}}
	\\
	\& 
	{\stackquiver{
	{\left\{\drawStringOrBand{0}{a,d,c}{3}\right\}}
	\\ {\scriptstyle\textcolor{teal!50!blue}{(0, 0, 1, 1, 1, 0, 0)}}}}
	\&\& 
	{\stackquiver{
	{\left\{\drawStringOrBand{0}{b,a,d}{3}\right\}}
	\\ {\scriptstyle\textcolor{teal!50!blue}{(0, 1, 1, 2, 0, 0, 0)}}}}
	\& 
	{\stackquiver{
	{\left\{\drawStringOrBand{1}{b,a}{2}, \drawStringOrBand{1}{d,c}{2}\right\}}
	\\ {\scriptstyle\textcolor{teal!50!blue}{(0, 1, 1, 1, 0, 0, 1)}}}}
	\& 
	{\stackquiver{
	{\left\{\drawStringOrBand{0}{d,c,b}{3}\right\}}
	\\ {\scriptstyle\textcolor{teal!50!blue}{(0, 0, 1, 1, 0, 1, 0)}}}}
	\&\& 
	{\stackquiver{
	{\left\{\drawStringOrBand{0}{c,b,a}{3}\right\}}
	\\ {\scriptstyle\textcolor{teal!50!blue}{(0, 1, 2, 1, 0, 0, 0)}}}}
	\\
	\&\&\&\& 
	{\stackquiver{
	{\left\{\drawStringOrBand{1}{b,a,d,c}{4}\right\}}
	\\ {\scriptstyle\textcolor{teal!50!blue}{(0, 0, 1, 1, 0, 0, 0)}}}}
	\arrow[no head, from=1-5, to=2-2]
	\arrow[no head, from=1-5, to=2-4]
	\arrow[no head, from=1-5, to=2-6]
	\arrow[no head, from=1-5, to=2-8]
	\arrow[no head, from=2-2, to=3-1]
	\arrow[no head, from=2-2, to=3-3]
	\arrow[no head, from=2-2, to=3-4]
	\arrow[no head, from=2-4, to=3-1]
	\arrow[no head, from=2-4, to=3-6]
	\arrow[no head, from=2-4, to=3-7]
	\arrow[no head, from=2-6, to=3-3]
	\arrow[no head, from=2-6, to=3-6]
	\arrow[no head, from=2-6, to=3-9]
	\arrow[no head, from=2-8, to=3-4]
	\arrow[no head, from=2-8, to=3-7]
	\arrow[no head, from=2-8, to=3-9]
	\arrow[squiggly, no head, from=3-1, to=4-1]
	\arrow[no head, from=3-1, to=4-2]
	\arrow[no head, from=3-1, to=4-4]
	\arrow[no head, from=3-3, to=4-2]
	\arrow[squiggly, no head, from=3-3, to=4-3]
	\arrow[no head, from=3-3, to=4-6]
	\arrow[no head, from=3-4, to=4-4]
	\arrow[no head, from=3-4, to=4-6]
	\arrow[no head, from=3-6, to=4-2]
	\arrow[no head, from=3-6, to=4-8]
	\arrow[no head, from=3-7, to=4-4]
	\arrow[squiggly, no head, from=3-7, to=4-7]
	\arrow[no head, from=3-7, to=4-8]
	\arrow[no head, from=3-9, to=4-6]
	\arrow[no head, from=3-9, to=4-8]
	\arrow[squiggly, no head, from=3-9, to=4-9]
	\arrow[no head, from=4-1, to=5-2]
	\arrow[no head, from=4-1, to=5-4]
	\arrow[squiggly, no head, from=4-2, to=5-2]
	\arrow[no head, from=4-2, to=5-5]
	\arrow[no head, from=4-3, to=5-2]
	\arrow[no head, from=4-3, to=5-6]
	\arrow[squiggly, no head, from=4-4, to=5-4]
	\arrow[no head, from=4-4, to=5-5]
	\arrow[no head, from=4-6, to=5-5]
	\arrow[squiggly, no head, from=4-6, to=5-6]
	\arrow[no head, from=4-7, to=5-4]
	\arrow[no head, from=4-7, to=5-8]
	\arrow[no head, from=4-8, to=5-5]
	\arrow[squiggly, no head, from=4-8, to=5-8]
	\arrow[no head, from=4-9, to=5-6]
	\arrow[no head, from=4-9, to=5-8]
	\arrow[no head, from=5-2, to=6-5]
	\arrow[no head, from=5-4, to=6-5]
	\arrow[squiggly, no head, from=5-5, to=6-5]
	\arrow[no head, from=5-6, to=6-5]
	\arrow[no head, from=5-8, to=6-5]
\end{tikzcd}
\end{equation*}
        \caption{Poset spanned by degenerations of strings and bands for dimension vector $(1,2,1)$. The squiggly lines represent the resolutions of reachings. The coloured vectors show certain entries of the $h$-vector.}
        \label{fig:deg_dim_121}
    \end{figure}
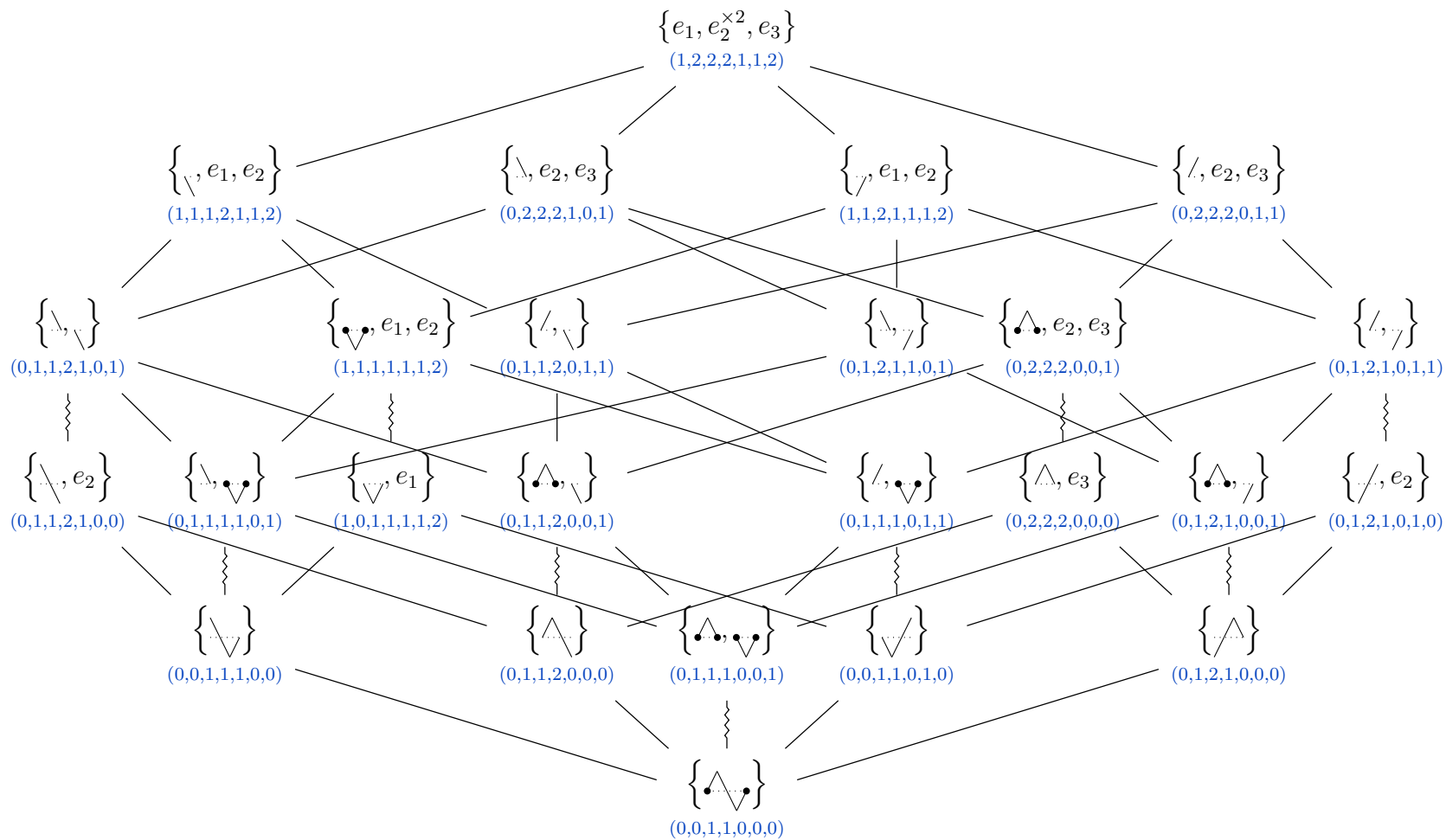
\end{landscape}

\begin{example}
\label[example]{ex:glued_Kronecker_dim242}
    We continue \Cref{ex:glued_Kronecker_dim121} with the dimension vector $(2,4,2)$. The poset of all degenerations becomes too big for discussion here, so we restrict to the mutli-sets only containing minimal bands and only consider degenerations arising from resolutions of reachings. The associated poset is given by \Cref{fig:ex_glued_Kronecker_band_deg_dim_242}. Note that the degeneration from \Cref{ex:resolution_produces_non_minimal_bands}, marked by the dashed line, is not a cover relation.
    \begin{figure}[H]
        \centering
        \begin{equation*}
\begin{tikzcd}
	& \left\{ \drawStringOrBand{1}{b,a}{2}^{\times 2}, \drawStringOrBand{1}{d,c}{2}^{\times 2} \right\} 
	\\
	&
 \left\{ \drawStringOrBand{1}{b,a,d,c}{4}, \drawStringOrBand{1}{b,a}{2}, \drawStringOrBand{1}{d,c}{2} \right\} 
 \\
 \left\{ \drawStringOrBand{1}{b,a,d,c,d,c}{6}, \drawStringOrBand{1}{b,a}{2}\right\} 
    &&
    \left\{ \drawStringOrBand{1}{b,a,b,a,d,c}{6}, \drawStringOrBand{1}{d,c}{2} \right\} 
    \\
 & \left\{ \drawStringOrBand{1}{b,a,b,a,d,c,d,c}{8} \right\} & 
 \\
 & \left\{ \drawStringOrBand{1}{b,a,d,c}{4}^{\times 2} \right\} 
	\arrow[no head, from=1-2, to=2-2]
	\arrow[no head, from=2-2, to=3-1]
	\arrow[no head, from=2-2, to=3-3]
	\arrow[no head, from=3-1, to=4-2]
	\arrow["\text{\Cref{ex:resolution_produces_non_minimal_bands}}"', dashed, no head, from=3-1, to=5-2]
	\arrow[no head, from=3-3, to=4-2]
	\arrow[no head, from=4-2, to=5-2]
\end{tikzcd}
\end{equation*}
        \caption{The poset of degenerations of multi-sets of bands arising from reachings for dimension vector $(2,4,2)$.}
        \label{fig:ex_glued_Kronecker_band_deg_dim_242}
    \end{figure}
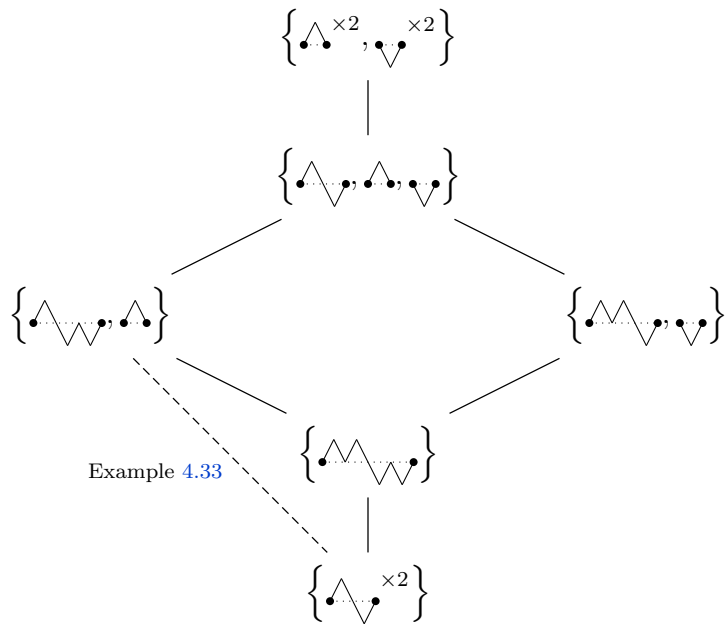
\end{example}

\begin{example}
\label[example]{ex:glued_Kronecker_dim363}
    We continue \Cref{ex:glued_Kronecker_dim121} and \Cref{ex:glued_Kronecker_dim242} with the dimension vector $(3,6,3)$. Again, we restrict to degenerations between the multi-sets of minimal bands arising from resolutions. The poset is shown in \Cref{fig:ex_glued_Kronecker_band_deg_dim_363}. Note that the top part, surrounded by dashed lines, is isomorphic to the poset of \Cref{ex:glued_Kronecker_dim242}.
    
    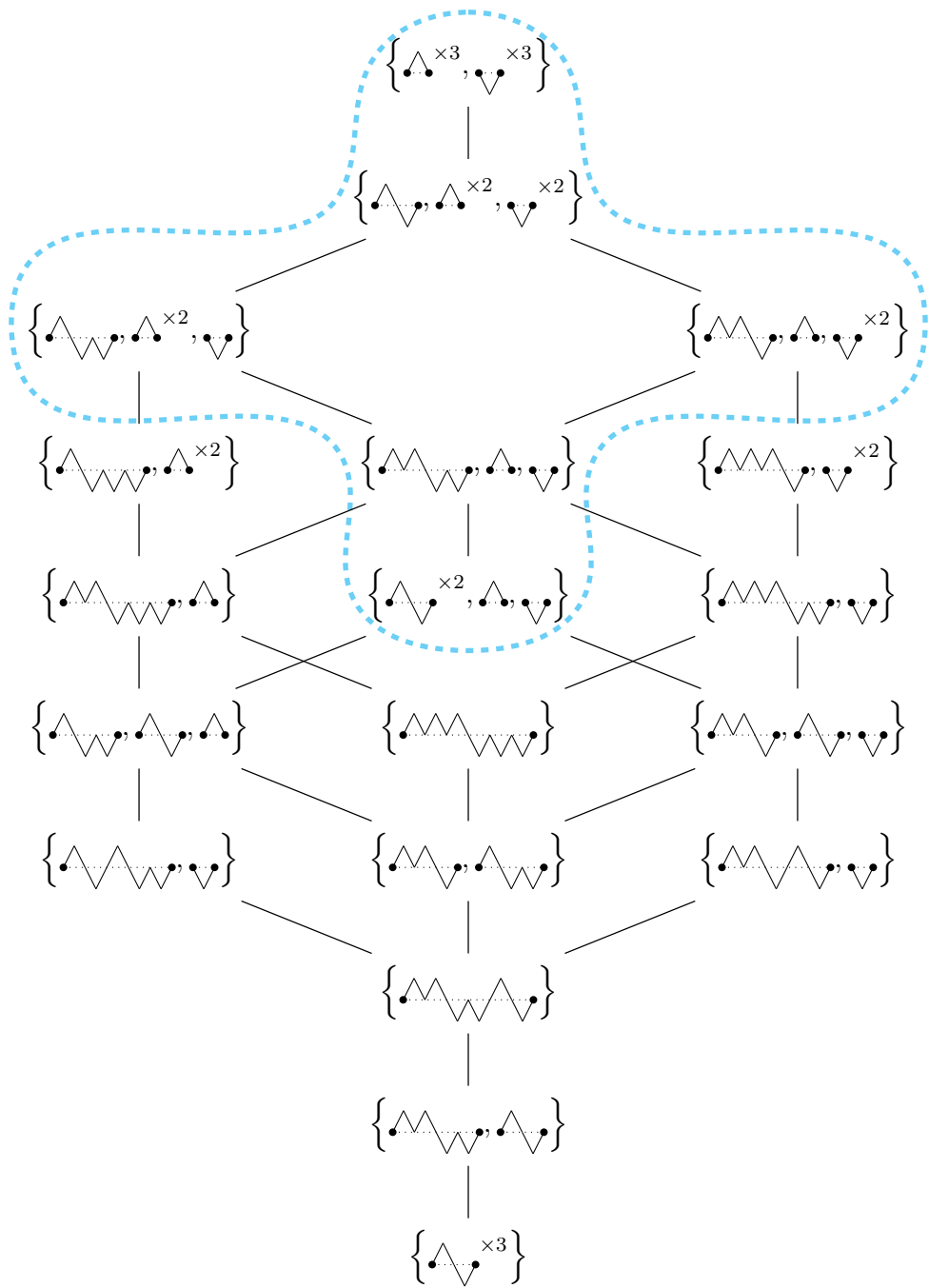
\begin{figure}[H]
        \centering
        \begin{equation*}
    \begin{tikzcd}[
execute at end picture={
  \begin{pgfonlayer}{background}
        \coordinate (a-l) at (A.west);
        \coordinate (a-t) at ($(A.north)+(0,0.2)$);
        \coordinate (a-r) at (A.east);
        \coordinate (a-b) at (A.south);
        \coordinate (b-l) at (B.west);
        \coordinate (b-r) at (B.east);
        \coordinate (c-l) at (C.west);
        \coordinate (d-r) at (D.east);
        \coordinate (e-l) at (E.west);
        \coordinate (e-r) at (E.east);
        \coordinate (f-l) at (F.west);
        \coordinate (f-r) at (F.east);
        \coordinate (f-b) at ($(F.south)+(0,-0.2)$);
        \draw[cyan!50!white, line width = 2 pt, dashed] (a-t) to[closed, curve through =
    { (a-r) (b-r) (d-r) (e-r) (f-r) (f-b) (f-l) (e-l) (c-l) (b-l) }] (a-l);
    \end{pgfonlayer}
}
]
	& |[alias=A]|{\left\{ \drawStringOrBand{1}{b,a}{2}^{\times 3}, \drawStringOrBand{1}{d,c}{2}^{\times 3} \right\}} 
	\\
	&
 |[alias=B]|{\left\{ \drawStringOrBand{1}{b,a,d,c}{4}, \drawStringOrBand{1}{b,a}{2}^{\times 2}, \drawStringOrBand{1}{d,c}{2}^{\times 2} \right\}}
 \\
 |[alias=C]|{\left\{ \drawStringOrBand{1}{b,a,d,c,d,c}{6}, \drawStringOrBand{1}{b,a}{2}^{\times 2}, \drawStringOrBand{1}{d,c}{2} \right\}}
    &&
    |[alias=D]|{\left\{ \drawStringOrBand{1}{b,a,b,a,d,c}{6}, \drawStringOrBand{1}{b,a}{2}, \drawStringOrBand{1}{d,c}{2}^{\times 2} \right\}}
    \\
	 \left\{ \drawStringOrBand{1}{b,a,d,c,d,c,d,c}{8}, \drawStringOrBand{1}{b,a}{2}^{\times 2} \right\} 
 &
  |[alias=E]|{\left\{ \drawStringOrBand{1}{b,a,b,a,d,c,d,c}{8}, \drawStringOrBand{1}{b,a}{2}, \drawStringOrBand{1}{d,c}{2} \right\}}
 & 
   \left\{ \drawStringOrBand{1}{b,a,b,a,b,a,d,c}{8}, \drawStringOrBand{1}{d,c}{2}^{\times 2} \right\} 
 \\
 \left\{ \drawStringOrBand{1}{b,a,b,a,d,c,d,c,d,c}{10}, \drawStringOrBand{1}{b,a}{2} \right\} 
 & 
  |[alias=F]|{\left\{ \drawStringOrBand{1}{b,a,d,c}{4}^{\times 2}, \drawStringOrBand{1}{b,a}{2}, \drawStringOrBand{1}{d,c}{2} \right\}}
    & 
 \left\{ \drawStringOrBand{1}{b,a,b,a,b,a,d,c,d,c}{10}, \drawStringOrBand{1}{d,c}{2} \right\} 
 \\
 \left\{ \drawStringOrBand{1}{b,a,d,c,d,c}{6}, \drawStringOrBand{1}{b,a,d,c}{4}, \drawStringOrBand{1}{b,a}{2} \right\} 
 &
 \left\{ \drawStringOrBand{1}{b,a,b,a,b,a,d,c,d,c,d,c}{12} \right\} 
 & 
  \left\{ \drawStringOrBand{1}{b,a,b,a,d,c}{6}, \drawStringOrBand{1}{b,a,d,c}{4}, \drawStringOrBand{1}{d,c}{2} \right\} 
 \\
	\left\{ \drawStringOrBand{1}{b,a,d,c,b,a,d,c,d,c}{10}, \drawStringOrBand{1}{d,c}{2} \right\}
    &
    \left\{ \drawStringOrBand{1}{b,a,b,a,d,c}{6}, \drawStringOrBand{1}{b,a,d,c,d,c}{6} \right\}
    &
     \left\{ \drawStringOrBand{1}{b,a,b,a,d,c,b,a,d,c}{10}, \drawStringOrBand{1}{d,c}{2} \right\}
    \\
	&
  \left\{ \drawStringOrBand{1}{b,a,b,a,d,c,d,c,b,a,d,c}{12} \right\}\\
	&  \left\{ \drawStringOrBand{1}{b,a,b,a,d,c,d,c}{8}, \drawStringOrBand{1}{b,a,d,c}{4} \right\} \\
	&  \left\{ \drawStringOrBand{1}{b,a,d,c}{4}^{\times 3} \right\} 
	\arrow[no head, from=1-2, to=2-2]
	\arrow[no head, from=2-2, to=3-1]
	\arrow[no head, from=2-2, to=3-3]
	\arrow[no head, from=3-1, to=4-1]
	\arrow[no head, from=3-1, to=4-2]
	\arrow[no head, from=3-3, to=4-2]
	\arrow[no head, from=3-3, to=4-3]
	\arrow[no head, from=4-1, to=5-1]
    \arrow[no head, from=4-2, to=5-1]
    \arrow[no head, from=4-2, to=5-2]
	\arrow[no head, from=4-2, to=5-3]
	\arrow[no head, from=4-3, to=5-3]
    \arrow[no head, from=5-1, to=6-1]
	\arrow[no head, from=5-1, to=6-2]
	\arrow[no head, from=5-2, to=6-1]
	\arrow[no head, from=5-2, to=6-3]
	\arrow[no head, from=5-3, to=6-2]
    \arrow[no head, from=5-3, to=6-3]
	\arrow[no head, from=6-1, to=7-1]
	\arrow[no head, from=6-1, to=7-2]
	\arrow[no head, from=6-2, to=7-2]
	\arrow[no head, from=6-3, to=7-2]
	\arrow[no head, from=6-3, to=7-3]
	\arrow[no head, from=7-1, to=8-2]
	\arrow[no head, from=7-2, to=8-2]
	\arrow[no head, from=7-3, to=8-2]
	\arrow[no head, from=8-2, to=9-2]
	\arrow[no head, from=9-2, to=10-2]
\end{tikzcd}
\end{equation*}
    
        \caption{The poset of degenerations of multi-sets of bands arising from reachings for the dimension vector $(3,6,3)$.}
        \label{fig:ex_glued_Kronecker_band_deg_dim_363}
    \end{figure}
\end{example}
    
\biblio{}

\section{Outlook}
\label{sect:outlook}

We leave with some questions and directions for possible future research. Already noted in \Cref{sect:degeneration}, the question if $\leqdeg$ coincides with the partial order induced by $h$ remains open. Closely related, we may ask if the orbits defined in \Cref{sect:families} for mutli-sets of strings and minimal bands, give a good stratification of $\mod(A,\dd)$, i.e. if $\sO_{\sD'} \cap \overline{\sO_\sD} \neq \varnothing$ implies that $\sO_{\sD'} \subseteq\overline{\sO_\sD}$. After \Cref{lem:deleting_arrow}, \Cref{thm::deg_and_intersection} and the examples that follow, it is natural to ask if every cover relation is induced by resolutions of reachings or deleting an arrow. This is true in all examples examined in this paper. Next, \Cref{rmk:families_string_alg,,ex:glued_Kronecker_dim242,,ex:glued_Kronecker_dim121} asks if it is possible to classify which resolutions of reachings or arrow deletions will give a cover relation. 

As noted in \Cref{rmk:h_vector_string_alg,,rmk:deletion_string_alg,,rmk:resolution_string_alg}, most results of this work can be easily generalised to string algebras. Since string and band modules can be defined in any basic algebra, we could even define string and band orbits for them and study their order.

\biblio{}

\section*{Acknowledgements}
This work started as a part of my master's thesis at the University of Bonn. I want to thank my thesis advisor Jan Schröer for the wonderful project proposal and his input and direction during the work on the thesis. Further, I thank René Marczinzik for excellent bibliography advice and Markus Kleinau for interesting discussions.

This work was continued during my PhD. I would like to thank my advisors Pierre-Guy Plamondon and Claire Amiot for their direction and insights.

This project has received funding from the European Union’s Horizon Europe research and innovation programme under the Marie Skłodowska-Curie grant agreement No 101126554.

\section*{Disclaimer}
Co-Funded by the European Union. Views and opinions expressed are however those of the author only and do not necessarily reflect those of the European Union. Neither the European Union nor the granting authority can be held responsible for them.
\noindent \includegraphics[height = 1cm]{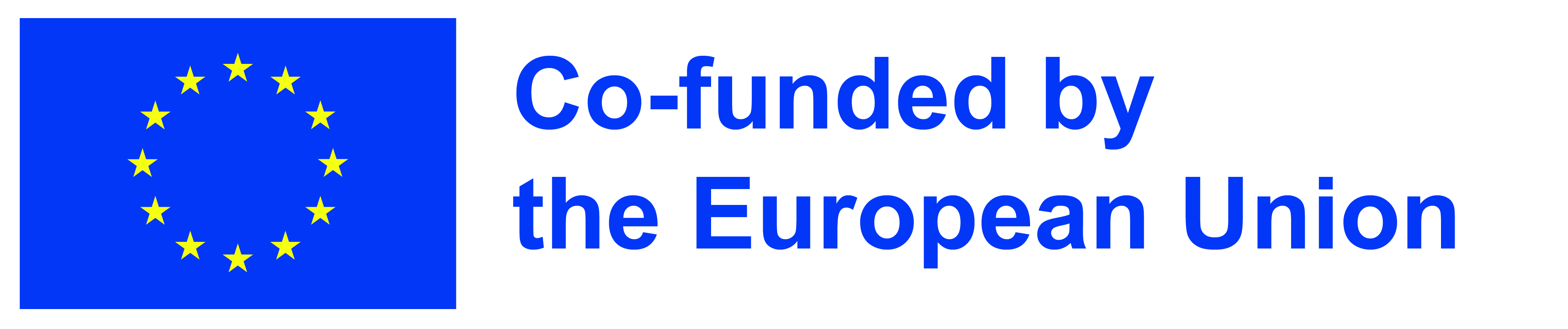}

\biblio

\bibliographystyle{amsalpha}
\bibliography{literature}

\end{document}